\documentclass[a4paper,12pt]{article}
\usepackage{amsmath, amsthm, amssymb}
\usepackage{graphicx,subfigure}
\usepackage[geometry]{ifsym}
\usepackage{anysize}
\usepackage{float}
\usepackage{caption}
\marginsize{2.5cm}{2.5cm}{2cm}{2cm}

\newtheorem{theorem}{Theorem}[section]
\newtheorem{lemma}[theorem]{Lemma}

\newtheorem{property}[theorem]{Property}
\newtheorem{corollary}[theorem]{Corollary}
\newtheorem{definition}[theorem]{Definition}

\begin{document}
\newcommand{\Eqref}[1]{(\ref{#1})}

\newcommand\scalemath[2]{\scalebox{#1}{\mbox{\ensuremath{\displaystyle #2}}}}

\def\CC{{\rm C\kern-.18cm\vrule width.6pt height 6pt depth-.2pt
\kern.18cm}}

\def\NN{{\mathop{{\rm I}\kern-.2em{\rm N}}\nolimits}}

\def\ZZ{{\mathbb{Z}}}

\def\PP{{\mathop{{\rm I}\kern-.2em{\rm P}}\nolimits}}

\def\RR{{\mathop{{\rm I}\kern-.2em{\rm R}}\nolimits}}

\title{Spaces of generalized splines over T-meshes}

\author{{Cesare Bracco}$^a$ and {Fabio Roman}$^a$
\\ \\
\small $^a$Department of Mathematics \lq\lq G. Peano\rq\rq - University of Turin\\ \small V. Carlo Alberto 10, Turin 10123, Italy\\
}\date{}
\maketitle


 

\begin{abstract} 

We consider a class of non-polynomial spline spaces over T-meshes, that is, of spaces locally spanned both by polynomial and by suitably-chosen non-polynomial functions, which we will refer to as generalized splines over T-meshes. For such spaces, we provide, under certain conditions, a dimension formula and a basis based on the notion of minimal determining set. We explicitly examine some relevant cases, which enjoy a noteworthy behaviour with respect to differentiation and integration; finally, we also study the approximation power of the just constructed spline spaces.
\end{abstract}



\section{Introduction} \label{sec:100}

\medskip
The T-splines, first introduced in \cite{sederberg03} and \cite{sederberg04}, represent a significant advancement in CAD and CAGD techniques, with relevant applications to differential problems, in particular in the framework of isogeometric analysis (see, e.g., \cite{cottrell} and \cite{bazilevs10}). The spline spaces over T-meshes are a closely related notion, first introduced by Deng et al. in \cite{deng06} and further studied by the same authors and several others (see, e.g., \cite{li06}, \cite{li11}, \cite{mourrain} and \cite{schum12}). The basic idea consists of considering spaces of spline functions which are polynomials of a certain degree in each of the cells of the T-mesh, which, unlike the classical tensor-product meshes, allows T-junctions, that is, vertices where only three edges meet.

\smallskip
\noindent 
This paper is devoted to the study of the {\sl Generalized spline spaces over T-meshes}, a class of non-polynomial spline spaces, which essentially generalize the concept of spline space over a T-mesh: roughly speaking, they are locally spanned by polynomials and some suitable non-polynomial functions. The relevance of this class of spline spaces and some of the basic concepts related to it have been only recently discussed in some international conferences. The same kind of non-polynomial functions have been recently used also to construct non-polynomial T-splines (see, e.g., \cite{bracco12b}), and non-polynomial hierarchical splines spaces (see \cite{manni14}). The goal of this work is to carry out a rigorous and deep study of the class of Generalized splines over T-meshes, including general results about the space dimension, the approximation power and noteworthy cases, which, as far as we know, are still missing in the literature. The introduction and the study of these non-polynomial spaces is justified by at least two reasons. First of all, the presence of non-polynomial functions allows to exactly reproduce certain shapes which can only be approximated by polynomial splines or NURBS (for example relevant curves like helices, cycloids, catenaries, or other transcendental curves). Moreover, as we will also point out in Section 4, choosing suitable non-polynomial functions also allows an easier computation of derivatives and integrals of certain surfaces with respect to using NURBS (see also \cite{manni11}, \cite{kvasov99}).

\smallskip
\noindent
Starting from some of the results obtained in \cite{cost05} about a noteworthy class of univariate non-polynomial spaces, we define generalized spline spaces over T-meshes and construct a local representation in the Bernstein-B\'ezier fashion for their elements. For the above spaces we first provide the construction of a basis and a dimension formula by using the properties of the local Bernstein-B\'ezier representation and by generalizing to the non-polynomial case the techniques proposed for the polynomial one in \cite{schum12}. 

\noindent
Then, we show and discuss some noteworthy cases of generalized spline spaces over a T-mesh, which have a good behaviour w.r.t. differentiation and integration. Such a feature is very useful for the applications of the considered spaces, especially in the isogeometric analysis.

\noindent
Moreover, we also study the approximation power of the just constructed spline spaces. In particular, we do it by constructing a quasi-interpolant based on some new local approximants, whose construction is not trivial. In fact, the results about the univariate non-polynomial Hermite interpolants given in \cite{cost05} cannot be directly extended to the bivariate case. On the other hand, also the bivariate averaged Taylor expansions used in \cite{schum12} cannot be simply adapted to the non-polynomial case we consider here. Therefore, we instead defined a new local Hermite interpolant belonging to the non-polynomial spline space, whose existence is proved by using certain assumptions made about the non-polynomial functions spanning the space, as carefully explained in Section 5. This approach allows us to get essentially the same approximation order as in the polynomial case. 

\smallskip
\noindent
The paper is organized as follows. Section \ref{sec:200} includes several preliminary arguments about the non-polynomial spaces we will use to define the new spline spaces, including some important properties about the derivatives of the basis functions and the basic concepts about T-meshes. Section \ref{sec:300} presents the new generalized spline spaces over T-meshes, and includes a detailed proof of the dimension formula and of the construction of the basis. Section \ref{sec:400} deals with the issue of suitable choices of the spaces in order to get a good behaviour of the spaces themselves w.r.t. differentiation and integration. Finally, Section \ref{sec:500} is devoted to the study of the approximation power of the constructed generalized spline space.

\section{Preliminaries} \label{sec:200}

The spaces we will consider are of the type
\begin{equation}
{\cal P}^n_{u,v}([a,b]):=\hbox{span}\langle 1,s,...,s^{n-3},u(s),v(s)\rangle, \qquad s\in[a,b], \quad 3\le n\in \NN, \label{p1ma1} 
\end{equation}
where $u,v\in C^n([a,b])$; for $n=2$ we set 
\begin{equation*}
{\cal P}^2_{u,v}([a,b]):=\hbox{span}\langle u(s),v(s)\rangle, \qquad s\in[a,b]. 
\end{equation*}
We assume that $\hbox{dim}\big({\cal P}^n_{u,v}([a,b])\big)=n$; moreover, in order to prove some of the properties we are about to present, we will sometimes require the following additional conditions on ${\cal P}^n_{u,v}([a,b])$ 
\begin{align}
\forall \psi\in {\cal P}^n_{u,v}([a,b]),\,\,&\hbox{if}\,\, \psi^{(n-2)}(s_1)=\psi^{(n-2)}(s_2)=0, \quad s_1,s_2\in[a,b],\,\, s_1\ne s_2\notag\\ 
&\hbox{then}\,\, \psi^{(n-2)}(s)=0, \quad s\in[a,b]; \label{p1ma2} 
\end{align}
\begin{align}
\forall \psi\in {\cal P}^n_{u,v}([a,b]),\,\,&\hbox{if}\,\, \psi^{(n-2)}(s_1)=\psi^{(n-1)}(s_1)=0, \quad s_1\in(a,b),\notag\\
&\hbox{then}\,\, \psi^{(n-2)}(s)=0, \quad s\in[a,b]. \label{p1ma3} 
\end{align}
In the following, we will explicitly mention when such conditions are needed.

\subsection{Normalized positive basis and its properties} \label{sec:210}

\medskip
In this subsection we consider a normalized positive basis for the space ${\cal P}^n_{u,v}([a,b])$. The procedure to obtain it and its fundamental properties are known and can be found in \cite{cost05}. Therefore here we will just recall the main results obtained in \cite{cost05}, omitting the proofs. We will instead prove Property \ref{p1ma13}, which will be crucial in order to obtain some results later in the paper. \\
In this Section we assume that the condition \eqref{p1ma2} holds.

\smallskip
The normalized positive basis can be constructed by using the following integral recurrence relation. By \eqref{p1ma2}, there exist unique elements $U_{0,1,n-1}$ and $U_{1,1,n-1}$ belonging to $\hbox{span}\langle u^{(n-2)},v^{(n-2)}\rangle$ satisfying
\begin{align}
U_{0,1,n-1}(a)=1,\quad U_{0,1,n-1}(b)=0,\notag\\
U_{1,1,n-1}(a)=0,\quad U_{1,1,n-1}(b)=1,\label{p1ma6} 
\end{align}
and
\begin{equation}
U_{0,1,n-1}(s),U_{1,1,n-1}(s)>0, \qquad s\in(a,b). \label{p1ma7} 
\end{equation}
Moreover, we define, for $k=2,...,n-1$ and $n\ge 3$ 
\begin{align}
&U_{0,k,n-1}(s)=1-V_{0,k-1,n-1}(s)\notag\\
&U_{i,k,n-1}(s)=V_{i-1,k-1,n-1}(s)-V_{i,k-1,n-1}(s),\qquad 1\le i\le k-1\notag\\
&U_{k,k,n-1}(s)=V_{k-1,k-1,n-1}(s), \label{p1ma16} 
\end{align}
where
\begin{equation}
V_{i,k,n-1}(s)=\int_a^s U_{i,k,n-1}/d_{i,k,n-1}dt, \label{p1ma8} 
\end{equation}
and
\begin{equation*}
d_{i,k,n-1}(s)=\int_a^b U_{i,k,n-1}dt, 
\end{equation*}
for $i=0,...,k$, $k=1,...,n-2$. Note that \eqref{p1ma6} and \eqref{p1ma7} hold also in the particular case $n=2$, and then $U_{0,1,1}$ and $U_{1,1,1}$ are a positive basis for ${\cal P}_{u,v}^2([a,b])$.
The following results can be proved about the just defined functions.

\begin{theorem}\label{p1ma10}
For $k=2,...,n-1$ and $n\ge 3$, the set of functions $\{U_{0,k,n-1},...,U_{k,k,n-1}\}$ is a basis for the space 
\begin{equation*}
\hbox{span}\langle 1,s,...,s^{k-2},u^{(n-k-1)}(s),v^{(n-k-1)}(s)\rangle. 
\end{equation*}
Moreover, it is a normalized positive basis, that is, satisfies the conditions $\sum_{i=0}^k U_{i,k,n-1}(s)=1$ and $U_{i,k,n-1}(s)>0$ for $s\in(a,b)$, $i=0,...,k$.
\end{theorem}

\begin{corollary}\label{p1ma11}
The set of functions $\{U_{0,n-1,n-1},...,U_{n-1,n-1,n-1}\}$ is a normalized positive basis for the space ${\cal P}_{u,v}^n([a,b])$, $n\ge 3$, $U_{i,n-1,n-1}=B_{i,n-1}$, where $\{B_{i,n-1}\}_{i=0}^{n-1}$ satisfy $\sum_{i=0}^{n-1} B_{i,n-1}(s)=1$ and $B_{i,n-1}(s)>0$ for $s\in(a,b)$, $i=0,...,n-1$. For $n=2$, the set $\{U_{0,1,1},U_{1,1,1}\}$ is a positive basis of ${\cal P}_{u,v}^2([a,b])$. 
\end{corollary}

\smallskip
\noindent
Since in the case $n=2$ we cannot, in general, guarantee the construction of a normalized positive basis, in the following we will assume $n\ge 3$. As a consequence of the results given in Sections 4 and 6 of \cite{cost05}, we get the following property. 

\begin{property}\label{p1ma12}
For $i=0,...,k$, $k=2,...,n-1$ and $n\ge 3$, we have
\begin{align*}
U_{i,k,n-1}^{(j)}(a)=0, \qquad j=0,...,i-1,\\
U_{i,k,n-1}^{(j)}(b)=0,\qquad j=0,...,k-i-1.
\end{align*}
In particular, if we consider $k=n-1$, we have
\begin{align*}
B_{i,n-1}^{(j)}(a)=0, \qquad j=0,...,i-1,\\
B_{i,n-1}^{(j)}(b)=0, \qquad j=0,...,n-i-2.
\end{align*}
\end{property}

\begin{property}\label{p1ma13}
For $k=2,...,n-1$ and $n\ge 3$, we have
\begin{align*}
U_{i,k,n-1}^{(i)}(a)\neq0, \qquad i=0,...,k-1,\\
U_{i,k,n-1}^{(k-i)}(b)\neq0, \qquad i=1,...,k. 
\end{align*}
In particular, if we consider $k=n-1$, we have
\begin{align}
B_{i,n-1}^{(i)}(a)\neq0, \qquad i=0,...,n-2,\label{p1ma14}\\
B_{i,n-1}^{(n-i-1)}(b)\neq0, \qquad i=1,...,n-1. \label{p1ma15}
\end{align}
\end{property}

\noindent
{\bf Proof.} First, let us prove \eqref{p1ma14} by induction. For $k=2$, \eqref{p1ma14} holds, since from \eqref{p1ma6}, \eqref{p1ma16} and \eqref{p1ma8} we get
\begin{align*}
&U_{0,2,n-1}(a)=1-V_{0,1,n-1}(a)=1-\int_{a}^a U_{0,1,n-1}(t)/d_{0,1,n-1} dt=1-0=1,\\
&U_{1,2,n-1}^{(1)}(a)=D[V_{0,1,n-1}(s)-V_{1,1,n-1}(s)]_{s=a}=\frac{U_{0,1,n-1}(a)}{d_{0,1,n-1}}-\frac{U_{1,1,n-1}(a)}{d_{1,1,n-1}}=\frac{1}{d_{0,1,n-1}}-0\neq 0.\\
\end{align*}
Now, if \eqref{p1ma14} holds for $k$, it must be true for $k+1$ as well, since we have
\begin{align*}
&U_{0,k+1,n-1}(a)=1-V_{0,k,n-1}(a)=1-\int_{a}^a U_{0,k,n-1}(t)/d_{0,k,n-1} dt=1-0=1, \\
&U_{i,k+1,n-1}^{(i)}=\frac{U_{i-1,k,n-1}^{(i-1)}(a)}{d_{i-1,k,n-1}}-\frac{U_{i,k,n-1}^{(i-1)}(a)}{d_{i,k,n-1}}=\frac{U_{i-1,k,n-1}^{(i-1)}(a)}{d_{i-1,k,n-1}}\neq 0,
\end{align*}
where we used \eqref{p1ma16}, \eqref{p1ma8}, Property \ref{p1ma12} and the induction hypothesis. Analogously we can prove \eqref{p1ma15}. \hfill$\square$

\smallskip
\noindent
Note that the above constructed basis is not only normalized positive, but it is also a Bernstein basis.

\subsection{Some definitions on T-meshes} \label{sec:220}

\medskip
We will now recall the definition of T-mesh and of some related objects, using the notations of \cite{schum12}. Note that the concept of T-mesh we will consider here may slightly differ from other ones in the literature, such as the more general used in \cite{daveiga}, which allows the presence not only of {\sl T-junctions}, but of {\sl L-junctions} and {\sl I-junctions} as well. 

\begin{definition}\label{p2ma1}
A T-mesh is a collection of axis-aligned rectangles $\Delta=\{R_i\}_{i=1}^N$ such that the domain $\Omega\equiv\cup_i R_i$ is connected and any pair of rectangles (which we will call {\sl cells}) $R_i,R_j\in \Delta$ intersect each other only at points on their edges.
\end{definition}

\vspace{-0.5cm}
\begin{figure}[ht] 
\centering
\includegraphics[scale=0.60]{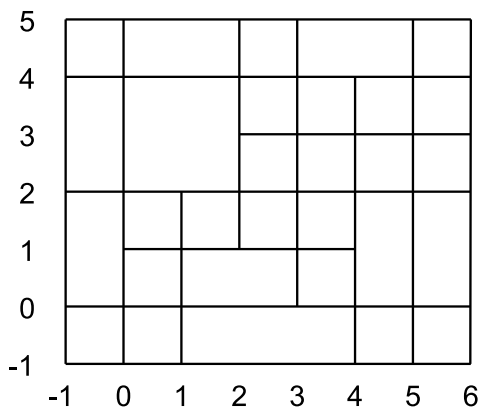}
\caption{{\sl An example of T-mesh where $\Omega=[-1,6]\times[-1,5]$}} \label{fig1}
\end{figure}

\smallskip
\noindent
Note that this definition does not imply that the domain $\Omega$ is rectangular and allows the presence of holes in it. Tensor-product meshes are a particular case of T-meshes. If a vertex $v$ of a cell belonging to $\Delta$ lies in the interior of an edge of another cell, then we call it a {\sl T-junction}.

\begin{definition}\label{p2ma2}
Given a T-mesh $\Delta$, a line segment $e=\langle w_1,w_2\rangle$ connecting the vertices $w_1$ and $w_2$ is called {\sl edge segment} if there are no vertices lying in its interior. Instead, if all the vertices lying in its interior are T-junctions and if it cannot be extended to a longer segment with the same property, then we call it a {\sl composite edge}. 
\end{definition}

In the following, we will consider T-meshes which are {\sl regular} and have no {\sl cycles}, in the sense of the following definitions (see \cite{schum12} for more details).

\begin{figure}[ht] 
\begin{minipage}[ht]{.43\textwidth}
\centering
\includegraphics[scale=0.42]{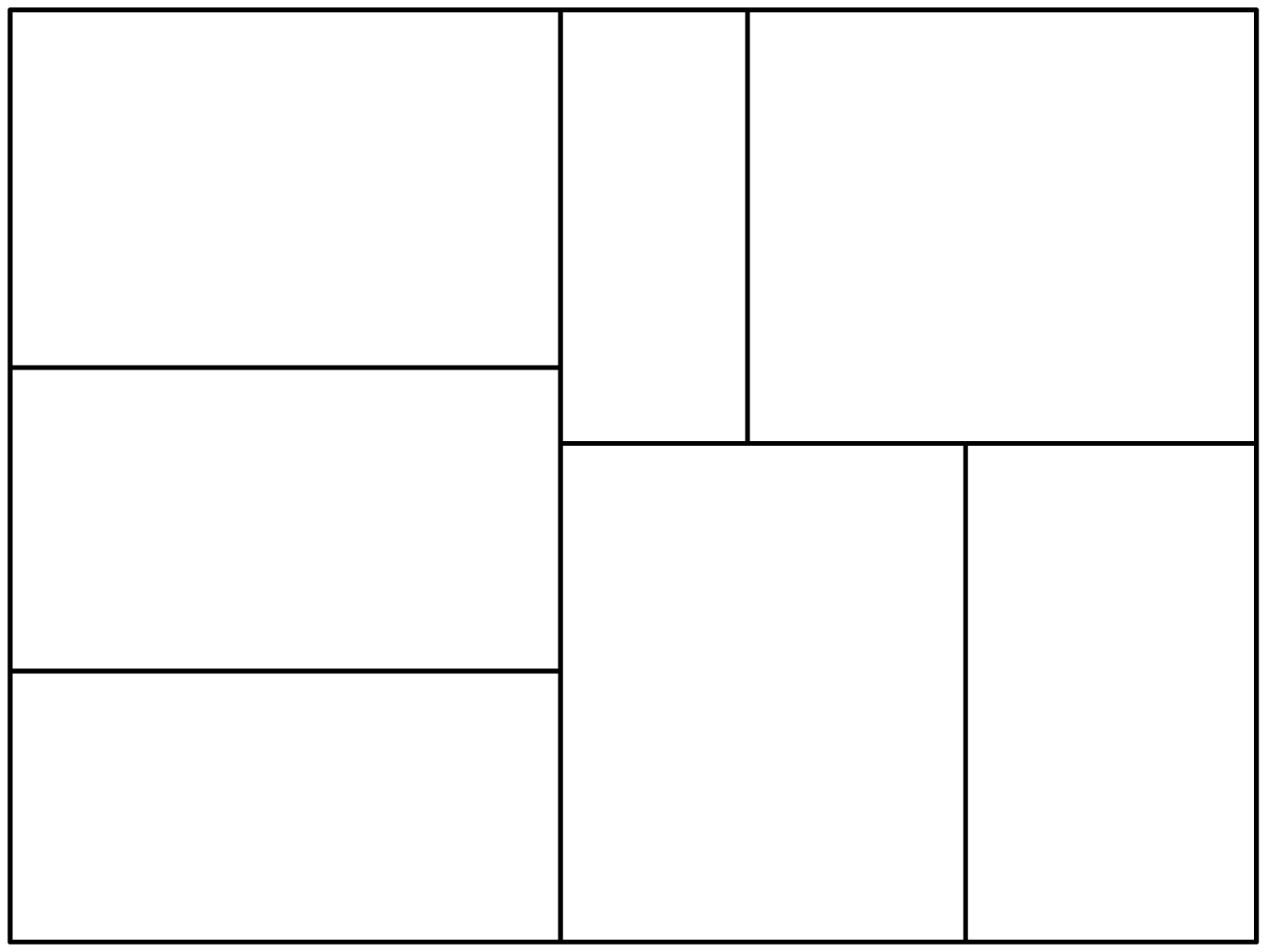}
\caption{\sl An example of regular T-mesh.} \label{fig2}
\end{minipage} 
\hspace{.1\textwidth}
\begin{minipage}[ht]{.43\textwidth}
\centering
\includegraphics[scale=0.42]{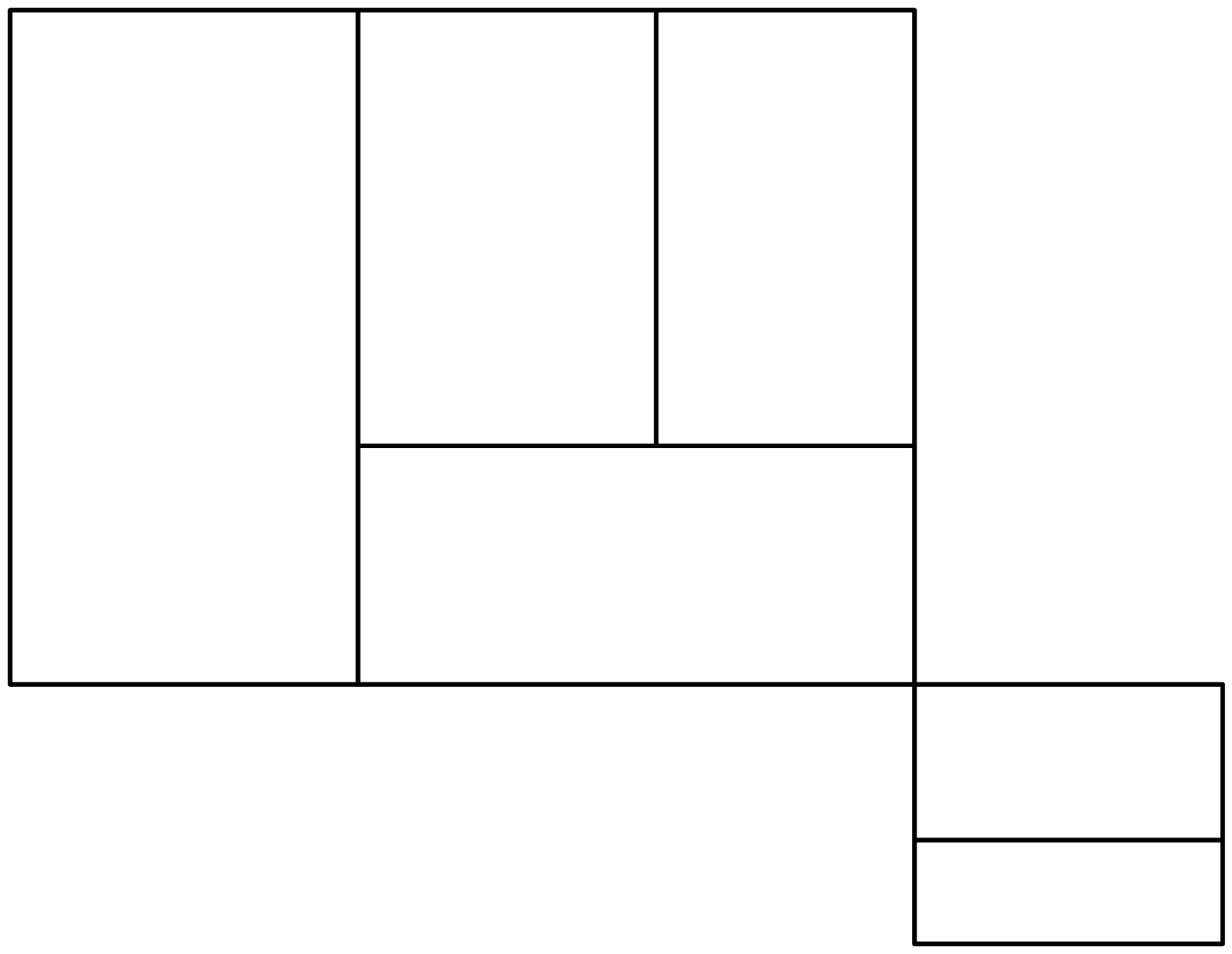}
\caption{\sl An example of non-regular T-mesh.} \label{fig3}
\end{minipage}
\end{figure}

\begin{definition}\label{p2ma3}
A T-mesh $\Delta$ is {\sl regular} if for each of its vertices $w$ the set of all rectangles containing $w$ has a connected interior. 
\end{definition}

\begin{definition}\label{p2ma4}
Let $w_1,...,w_n$ be a collection of T-junctions in a T-mesh $\Delta$ such that $w_i$ lies in the interior of a composite edge having one of its endpoints at $w_{i+1}$ (we assume $w_{n+1}=w_1$). Then $w_1,...,w_n$ are said to form a {\sl cycle}.
\end{definition}

\vspace{-0.5cm}
\begin{figure}[h] 
\centering
\includegraphics[scale=0.60]{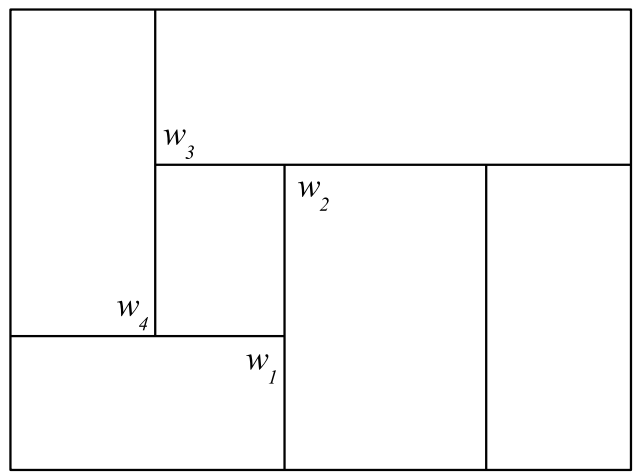}
\caption{{\sl The sequence $w_1,w_2,w_3,w_4$ is a cycle.}} \label{fig4}
\end{figure}

\section{Spaces of generalized splines on T-meshes} \label{sec:300}

\medskip
In this Section, we define the spaces of generalized splines over T-meshes, and we study their dimension by constructing a basis. The results obtained can be considered a generalization to non-polynomial splines spaces over T-meshes of the ones proved in \cite{schum12} for the basic polynomial case. 

\subsection{Basics} \label{sec:310}

\medskip
Let $\Delta$ be a regular T-mesh without cycles, and let $0\le r_1< n_1-1$, $0\le r_2< n_2-1$, where $r_1,r_2,n_1,n_2$ are integers and $n_1,n_2\ge2$. Later on, we will also use the notation ${\bf r}=(r_1,r_2)$ and ${\bf n}=(n_1,n_2)$. We define the {\sl space of generalized splines over the T-mesh $\Delta$ of bi-order ${\bf n}$ and smoothness ${\bf r}$} as
\begin{equation}
GS_{\bf u,v}^{\bf n,r}(\Delta)=\{p(s,t)\in C^{\bf r}(\Omega):\, p\vert_{R}\in {\cal P}^{\bf n}_{u,v}(R)\,\, \forall R\in \Delta\},\label{a2ma1}
\end{equation}
where $\Omega=\cup_{R\in \Delta} R$, $C^{\bf r}(\Omega)$ denotes the space of functions $p$ such that their derivatives $D^i_sD^j_tp$ are continuous for all $0\le i\le r_1$ and $0\le j\le r_2$, and the space ${\cal P}^{\bf n}_{\bf u,v}$ is defined as
\begin{equation}
{\cal P}^{\bf n}_{\bf u,v}(R) = \hbox{span} \langle g_1(s) g_2(t): \,
g_1 \in {\cal P}^{n_1}_{u_1,v_1}([a_R,b_R]), \,
g_2 \in {\cal P}^{n_2}_{u_2,v_2}([c_R,d_R]) \rangle,\label{a2ma2} 
\end{equation}
with $R=[a_R,b_R]\times[c_R,d_R]$ and $u_1,v_1\in C^{n_1}([a_R,b_R])$, $u_2,v_2\in C^{n_2}([c_R,d_R])$ such that $\hbox{dim}\big({\cal P}^{n_1}_{u_1,v_1}([a_R,b_R])\big)=n_1$ and \linebreak $\hbox{dim}\big({\cal P}^{n_2}_{u_2,v_2}([c_R,d_R])\big)=n_2$, and satisfying both \eqref{p1ma2} and \eqref{p1ma3}. In other words, $GS_{\bf u,v}^{\bf n,r}(\Delta)$ is a space of spline functions which, restricted to each cell $R$, are products of functions belonging to spaces of type \eqref{p1ma1}.

\bigskip
We introduce now on each cell $R$ a Bernstein-B\'ezier representation for the elements of $GS_{\bf u,v}^{\bf n,r}(\Delta)$ based on the Bernstein basis of ${\cal P}^{n_1}_{u_1,v_1}([a_R,b_R])$ and ${\cal P}^{n_2}_{u_2,v_2}([c_R,d_R])$ constructed in Theorem \ref{p1ma11}; therefore, we need to assume that \eqref{p1ma2} is satisfied both by ${\cal P}^{n_1}_{u_1,v_1}([a_R,b_R])$ and ${\cal P}^{n_2}_{u_2,v_2}([c_R,d_R])$. Let us denote by $\{B_{i,n_1-1}^R\}_{i=0}^{n_1-1}$ and $\{B_{i,n_2-1}^R\}_{i=0}^{n_2-1}$ the Bernstein basis of, respectively, ${\cal P}^{n_1}_{u_1,v_1}([a_R,b_R])$ and ${\cal P}^{n_2}_{u_2,v_2}([c_R,d_R])$, to stress the dependence of the basis on the coordinates $a_R,b_R,c_R,d_R$ of the vertices of the cell $R$. For any $p\in GS_{\bf u,v}^{\bf n,r}(\Delta)$, we can then give on the cell $R$ the following representation 
\begin{equation}
p\vert_R(s,t)=\sum_{i=0}^{n_1-1}\sum_{j=0}^{n_2-1}c_{ij}^R B_{i,n_1-1}^{R}(s)B_{j,n_2-1}^{R}(t),\label{a2ma3}
\end{equation}
where $c_{ij}^R\in \RR$ are suitable coefficients.
Let us define the set of {\sl domain points associated to R} 
\begin{equation*}
{\cal D}_{{\bf n},R}=\{\xi_{ij}^R\}_{i=0,j=0}^{n_1-1,n_2-1},
\end{equation*}
with
\begin{equation*}
\xi_{ij}^R=\Big(\frac{(n_1-1-i)a_R+ib_R}{n_1-1},\frac{(n_2-1-j)c_R+jd_R}{n_2-1}\Big), \qquad i=0,...,n_1-1,\,j=0,...,n_2-1.
\end{equation*}
We can then define the {\sl set of domain points} for a given T-mesh $\Delta$ as
\begin{equation*}
{\cal D}_{{\bf n},\Delta}=\bigcup_{R\in\Delta} {\cal D}_{{\bf n},R},
\end{equation*} 
where we assume that multiple appearances of the same point are allowed.
If we set
\begin{equation*}
B^R_\xi(s,t)=B_{i,n_1-1}^{R}(s)B_{j,n_2-1}^{R}(t), \quad\hbox{where} \,\, \xi_{ij}^R=\xi,  
\end{equation*}
then, for each $R\in\Delta$, we can re-write \eqref{a2ma3} in the more compact form
\begin{equation*}
p\vert_R(s,t)=\sum_{\xi\in{\cal D}_{{\bf n},R}}c_{\xi}^R B_{\xi}^{R}(s,t),
\end{equation*}
which we call {\sl Bernstein-B\'ezier form}; we refer to the $c_{\xi}^R$ as the {\sl B-coefficients}. It is then clear that any element of the space $GS_{\bf u,v}^{\bf n,r}(\Delta)$ is completely determined by a set of B-coefficients $\{c_\xi\}_{\xi\in{\cal D}_{{\bf n},\Delta}}$. Of course, not every choice of the B-coefficients corresponds to an element in the spline space, since smoothness conditions must be satisfied.

\subsection{Smoothness conditions} \label{sec:320}

\medskip
In order to study the consequences of the smoothness conditions required for $GS_{\bf u,v}^{\bf n,r}(\Delta)$ on the determination of the B-coefficients of an element of the space, first we need to recall some more concepts about domain points. 

\smallskip
Let $R=[a_R,b_R]\times[c_R,d_R]\in\Delta$, $w=(a_R,c_R)$, and ${\bf \mu}=(\mu_1,\mu_2)$ with $\mu_1\le n_1-1$ and $\mu_2\le n_2-1$. We call the set $ {\cal D}_{\bf \mu}^R(w)=\{\xi_{ij}\}_{i=0,j=0}^{\mu_1,\mu_2}$ the {\sl disk of size ${\bf \mu}$ around $w$}. The disks around the other vertices of $R$ can be defined analogously. Moreover, we say that the points $\xi_{ij}^R$ with $0\le i\le \nu$ lie within a distance $\nu$ from the edge $e=\{a_R\}\times[c_R,d_R]$ and we use the notation $d(\xi_{ij}^R)\le \nu$. Analogous notations hold for the other edges of $R$.

\noindent
Moreover, we can define the set of domain points
\begin{equation*}
{\cal D}_{\bf \mu}(w)=\bigcup_{R\in\Delta_w} {\cal D}_{\bf \mu}^R, 
\end{equation*}
where $\Delta_w\subset\Delta$ contains only the cells having $w$ as one of their vertices and multiple appearances of a point are allowed in the union. Given a composite edge $e$, an edge $\tilde e$ lying on $e$ and a vertex $w$ of $\tilde e$, if $d(w,\tilde e)\le \nu$, then we write that $d(w,e)\le \nu$ as well. 

\smallskip
\noindent
The following lemma is a key step to be able to understand the influence of the smoothness conditions and to get a dimension formula for the space.

\begin{lemma}\label{a2ma5}
Let $p\in GS_{\bf u,v}^{\bf n,r}(\Delta)$ and let $w$ be a vertex of $\Delta$. Let us consider two cells $R$ and $\tilde R$ with vertices (in counter-clockwise order) $w,w_2,w_3,w_4$ and $w,w_5,w_6,w_7$, respectively. If the coefficients $c_\xi$, $\xi \in {\cal D}_{\bf r}^R(w)$ are given, then the coefficients $c_\eta$, $\eta \in {\cal D}_{\bf r}^{\tilde R}(w)$ are uniquely determined by the smoothness conditions at $w$.
\end{lemma}

\smallskip
\noindent
{\bf Proof.} We prove the lemma for the case where $w$ is the upper-right corner of $R=[a_R,b_R]\times[c_R,d_R]$ and the lower-left corner of $\tilde R=[a_{\tilde R},b_{\tilde R}]\times[c_{\tilde R},d_{\tilde R}]$, that is, $w=(b_R,d_R)=(a_{\tilde R},c_{\tilde R})$. First, let us consider the partial derivatives of $p\vert_{\tilde R}$ 
\begin{equation*}
D_s^{h}D_t^{k}p\vert_{\tilde R}(a_{\tilde R},c_{\tilde R})=\sum_{i=0}^{n_1-1}\sum_{j=0}^{n_2-1}c_{ij}^{\tilde R} D_s^{h}B_{i,n_1-1}^{{\tilde R}}(a_{\tilde R})D_t^{k}B_{j,n_2-1}^{{\tilde R}}(c_{\tilde R}), \qquad 0\le h\le r_1,\,0\le k\le r_2.
\end{equation*}
Since by Corollary \ref{p1ma11} $B_{i,n_1-1}^{{\tilde R}}(s)=U_{i,n_1-1,n_1-1}(s)$ and $B_{j,n_2-1}^{{\tilde R}}(t)=U_{j,n_2-1,n_2-1}(t)$, using Property \ref{p1ma12} gives that 
\begin{align*}
&D_s^{h}B_{i,n_1-1}^{{\tilde R}}(a_{\tilde R})=0, \qquad h< i\le n_1-1,\\
&D_t^{k}B_{j,n_2-1}^{{\tilde R}}(c_{\tilde R})=0, \qquad k< j\le n_2-1.
\end{align*}
Therefore,
\begin{equation*}
D_s^{h}D_t^{k}p\vert_{\tilde R}(a_{\tilde R},c_{\tilde R})=\sum_{i=0}^{h}\sum_{j=0}^{k}c_{ij}^{\tilde R} D_s^{h}B_{i,n_1-1}^{{\tilde R}}(a_{\tilde R})D_t^{k}B_{j,n_2-1}^{{\tilde R}}(c_{\tilde R}). \label{a2ma6}
\end{equation*}
Now let us compute the partial derivatives of $p\vert_R$ 
\begin{equation*}
D_s^{h}D_t^{k}p\vert_R(b_R,d_R)=\sum_{i=0}^{n_1-1}\sum_{j=0}^{n_2-1}c_{ij}^R D_s^{h}B_{i,n_1-1}^{R}(b_R)D_t^{k}B_{j,n_2-1}^{R}(d_R), \qquad 0\le h\le r_1,\,0\le k\le r_2.
\end{equation*}
Since by Corollary \ref{p1ma11} $B_{i,n_1-1}^{R}(s)=U_{i,n_1-1,n_1-1}(s)$ and $B_{j,n_2-1}^{R}(t)=U_{j,n_2-1,n_2-1}(t)$, using Property \ref{p1ma12} gives that 
\begin{align*}
&D_s^{h}B_{i,n_1-1}^{R}(b_R)=0, \qquad 0\le i<n_1-1-h,\\
&D_t^{k}B_{j,n_2-1}^{R}(d_R)=0, \qquad 0\le j<n_2-1-k.
\end{align*}
Therefore,
\begin{equation*}
D_s^{h}D_t^{k}p\vert_R(b_R,d_R)=\sum_{i=n_1-1-h}^{n_1-1}\sum_{j=n_2-1-k}^{n_2-1}c_{ij}^R D_s^{h}B_{i,n_1-1}^{R}(b_R)D_t^{k}B_{j,n_2-1}^{R}(d_R). 
\end{equation*}
Requiring the $C^{\bf r}$ smoothness at $w$ is then equivalent to the linear system composed of the equations
\begin{equation}
\sum_{i=0}^{h}\sum_{j=0}^{k}c_{ij}^{\tilde R} D_s^{h}B_{i,n_1-1}^{{\tilde R}}(a_{\tilde R})D_t^{k}B_{j,n_2-1}^{{\tilde R}}(c_{\tilde R})=\sum_{i=n_1-1-h}^{n_1-1}\sum_{j=n_2-1-k}^{n_2-1}c_{ij}^R D_s^{h}B_{i,n_1-1}^{R}(b_R)D_t^{k}B_{j,n_2-1}^{R}(d_R),\label{a2ma8}
\end{equation}
for $h=0,...,r_1$, $k=0,...,r_2$.

\noindent
Note that in this case we have $\{c_{\xi}\}_{\xi \in {\cal D}_{\bf r}^R(w)}=\{c_{ij}^R\}_{i=n_1-1-r_1,j=n_2-1-r_2}^{n_1-1,n_2-1}$, that is, the $(r_1+1)\times(r_2+1)$ B-coefficients associated to $R$ given by hypothesis are exactly the ones on the right-hand of equations \eqref{a2ma8}. Analogously, $\{c_{\eta}\}_{\eta \in {\cal D}_{\bf r}^{\tilde R}(w)}=\{c_{ij}^{\tilde R}\}_{i=0,j=0}^{r_1,r_2}$, which means that the $(r_1+1)\times(r_2+1)$ B-coefficients associated to $\tilde R$ are the unknowns of the system \eqref{a2ma8}. It is easy to observe that, if we organize the equations according to the order of the derivatives, the matrix of the system is lower triangular. Moreover, the entries on the diagonal of the matrix, that is,
\begin{equation*}
D_s^{h}B_{h,n_1-1}^{\tilde R}(a_{\tilde R})D_t^{k}B_{k,n_2-1}^{\tilde R}(c_{\tilde R}), \qquad h=0,...,r_1,\,k=0,...,r_2,
\end{equation*} 
are not zero because of Property \ref{p1ma13}. \hfill$\square$

\medskip
After having studied the influence of smoothness around a vertex, we now study the situation around edges. Given a composite edge $e$, we will use the following notation 
\begin{equation*}
r_e=\begin{cases} r_1,&\text{if $e$ is vertical,}\\ r_2,&\text{if $e$ is horizontal,} \end{cases} 
\end{equation*}
\begin{equation*}
D_e=\begin{cases} D_s,&\text{if $e$ is vertical,}\\ D_t,&\text{if $e$ is horizontal,} \end{cases}
\end{equation*}
\begin{equation*}
n_e=\begin{cases} n_2,&\text{if $e$ is vertical,}\\ n_1,&\text{if $e$ is horizontal,} \end{cases} 
\end{equation*}
\begin{equation*}
u_e=\begin{cases} u_2,&\text{if $e$ is vertical,}\\ u_1,&\text{if $e$ is horizontal,} \end{cases} 
\end{equation*}
\begin{equation*}
v_e=\begin{cases} v_2,&\text{if $e$ is vertical,}\\ v_1,&\text{if $e$ is horizontal.} \end{cases} 
\end{equation*}
Moreover, to get the following results we will assume that $u_1,v_1$ are such that for each horizontal edge segment $f=[a_f,b_f]\times \{c_f\}$
\begin{equation}
\dim {\cal P}^{n_1}_{u_1,v_1}([a_f,b_f])=n_1, \label{p23ot1} 
\end{equation}
and that $u_2,v_2$ are such that for each vertical edge segment $f=\{a_f\}\times [c_f,d_f]$
\begin{equation}
\dim {\cal P}^{n_2}_{u_2,v_2}([c_f,d_f])=n_2, \label{p23ot2}. 
\end{equation}

\begin{lemma}\label{a7ma7}
Let $e$ be a composite edge of $\Delta$ with endpoints $w_1$ and $w_5$, and let $p\in GS_{\bf u,v}^{\bf n,r}(\Delta)$. For any $0\le \nu\le r_e$, $D_e^{\nu}p\vert_e$ is a univariate function belonging to the space ${\cal P}^{n_e}_{u_e,v_e}([w_{1,e},w_{5,e}])$, where $w_{1,e},w_{5,e}\in \RR$ are the abscissas/ordinates of $w_1,w_5$. 
\end{lemma}

\smallskip
\noindent
{\bf Proof.} The Lemma can be trivially proved with the same arguments used in \cite{schum12} for the polynomial case, thanks to the assumptions \eqref{p23ot1} and \eqref{p23ot2}. \hfill$\square$

\smallskip
Let us now consider a cell $R_e$ with vertices $w_1,w_2,w_3,w_4$ and another cell $\tilde R_e$ with vertices $w_5,w_6,w_7,w_8$. Moreover we assume that $w_2$ and $w_8$ lie on $e$ as well (the other cases are analogous). Let us use the notation
\begin{equation}
{\cal M}_e=\begin{cases} \Big\{\xi_{ij}^{R_e}\Big\}_{i=0,j=r_2+1}^{r_1,n_2-r_2-2},&\text{if $e$ is vertical,}\\\Big\{\xi_{ij}^{R_e}\Big\}_{i=r_1+1,j=0}^{n_1-r_1-2,r_2},&\text{if $e$ is horizontal.}\end{cases}\label{a8ma1}
\end{equation} 
In other words, the set ${\cal M}_e$ includes all the domain points $\xi$ lying outside the disks ${\cal D}_{\bf r}^{R_e}(w_1)$ and
${\cal D}_{\bf r}^{R_e}(w_2)$ and satisfying $d(\xi,e)\le r_e$.

\begin{lemma}\label{a8ma2}
Let $e$ be an edge of the T-mesh $\Delta$ with endpoints $w_1$ and $w_5$, and let us assume that ${\cal P}^{n_e}_{u_e,v_e}([w_{1,e},w_{5,e}])$ satisfies \eqref{p1ma2}. If the B-coefficients of a spline $p\in GS_{\bf u,v}^{\bf n,r}(\Delta)$ corresponding to the domain points belonging to the set
\begin{equation*}
\tilde{\cal M}_e={\cal D}_{\bf r}^{R^e}(w_1)\cup {\cal D}_{\bf r}^{\tilde R^e}(w_5)\cup {\cal M}_e 
\end{equation*}
are given, then the coefficient of $p$ associated to domain points $\xi$ such that $d(\xi,e)\le r_e$ are uniquely determined.
\end{lemma}

\begin{figure}[h] 
\centering
\includegraphics[scale=0.60]{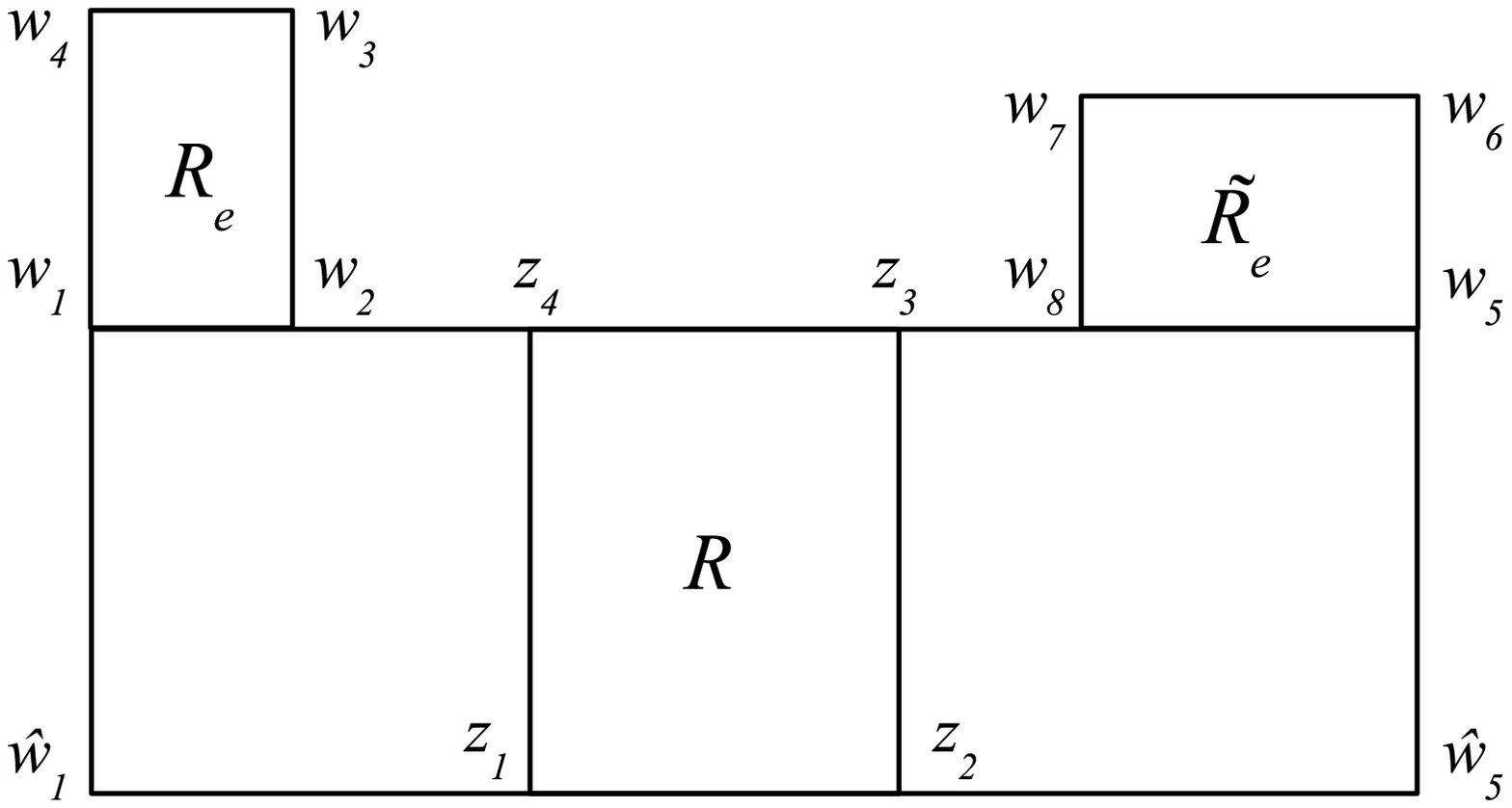}
\caption{{\sl the cells considered in the proof of Lemma \ref{a8ma2} (see \cite{schum12} for the analogous situation examined in the polynomial case).}} \label{fig5}
\end{figure}

\smallskip
\noindent
{\bf Proof.} We will suppose that $e$ is horizontal (the case where $e$ is vertical is analogous). We consider a cell $R$ with an edge lying on $e$: let us assume, for instance that $R$ has vertices $z_1,z_2,z_3,z_4$ and that $z_3$ and $z_4$ lie on $e$, like in Figure \ref{fig5} (see \cite{schum12} for the corresponding polynomial case). We will show that the B-coefficients corresponding to the domain points $\xi$ belonging to ${\cal D}_{{\bf n},R}$ and such that $d(\xi,e)\le r_2$ are uniquely determined. 

\noindent
Let $p\in GS_{\bf u,v}^{\bf n,r}(\Delta)$. First of all, since the B-coefficients corresponding to the domain points in $\tilde{\cal M}_e$ are given, we can compute the derivatives
\begin{equation*}
\Big\{D^i_sD^j_t p(w_1)\Big\}_{i=0,j=0}^{n_1-r_1-2,r_2}, \qquad \Big\{D^i_sD^j_t p(w_5)\Big\}_{i=0,j=0}^{r_1,r_2}.
\end{equation*}
In fact, by Property \ref{p1ma12}, the computation of these derivatives involves just the B-coefficients contained in $\tilde{\cal M}_e$. Now, let us consider the smallest rectangle $\tilde R$ containing $R$ and $e$ (the rectangle with vertices $\hat w_1,\hat w_5,w_5,w_1$ in Figure \ref{fig5}). Note that $\tilde R$ does not necessarily belong to the T-mesh $\Delta$; however, since we assumed that ${\cal P}^{n_e}_{u_e,v_e}([w_{1,e},w_{5,e}])$ satisfies \eqref{p1ma2}, we can temporarily assume it does and consider the corresponding B-coefficients. In order to obtain the B-coefficients $c^{\tilde R}_{ij}$, $i=0,...,n_1-1$, $j=0,...,r_2$, associated to $\tilde R$, that is, the ones corresponding to domain points $\xi$ in $\tilde R$ s.t. $d(\xi,e)\le r_e$, it is sufficient to solve the linear systems 
\begin{align}
&D^i_sD^j_tp\vert_{\tilde R}(w_1)=D^i_sD^j_tp(w_1), \qquad i=0,...,n_1-r_1-2\,\quad j=0,...,r_2,\notag\\
&D^i_sD^j_tp\vert_{\tilde R}(w_5)=D^i_sD^j_tp(w_5), \qquad i=0,...,r_1\,\quad j=0,...,r_2, \label{a13ma1}
\end{align}
where on the right side of the equality we have the already computed derivatives. Note that, by Properties \ref{p1ma12} and \ref{p1ma13}, the matrices associated to these systems are both triangular with non-zero elements on the diagonal. From Lemma \ref{a7ma7} we know that $p\vert_e\in {\cal P}^{n_e}_{u_e,v_e}([w_{1,e},w_{5,e}])$; in particular, this allows us to state that 
\begin{align}
&D^i_sD^j_tp\vert_{R}(w)=D^i_sD^j_tp\vert_{\tilde R}(w), \qquad i=0,...,n_1-r_1-2\,\quad j=0,...,r_2, \notag\\
&D^i_sD^j_tp\vert_{R}(\hat w)=D^i_sD^j_tp\vert_{\tilde R}(\hat w), \qquad i=0,...,r_1\,\quad j=0,...,r_2, \label{a13ma2}
\end{align}
for any two points $w$ and $\hat w$ lying on $e$ between $z_4$ and $z_3$. Let us choose $w=z_4$ and $\hat w=z_3$: we observe that the derivatives on the right side of equations \eqref{a13ma2} can be computed because we have enough B-coefficients of $p\vert_{\tilde R}$ (previously obtained solving \eqref{a13ma1}). Then, we get from \eqref{a13ma2} two linear systems where the unknowns are the B-coefficients $c^R_{ij}$, $i=0,...,n_1-1$, $j=0,...,r_2$, associated to the cell $R$, that is, the ones corresponding to domain points $\xi$ in $R$ s.t. $d(\xi,e)\le r_e$. These systems can be solved, since the associated matrices are triangular and with non-zero elements on the diagonals (again because of Properties \ref{p1ma12} and \ref{p1ma13}). Since the same procedure can be repeated for any cell with one edge lying on $e$, the Lemma is proved. \hfill$\square$

\subsection{Basis and dimension formula} \label{sec:330}

\medskip
We will prove the construction of the basis and a dimension formula for $GS_{\bf u,v}^{\bf n,r}(\Delta)$ in the case $n_1-1\ge 2r_1+1$ and $n_2-1\ge 2r_2+1$. One of the main tools which we will need to prove the dimension formula is the concept of {\sl minimal determining set}, introduced in \cite{schum12} for the polynomial case.

\begin{definition}\label{a7ma1}
Let ${\cal M}\subset {\cal D}_{{\bf n}, \Delta}$. ${\cal M}$ is a {\sl determining set} for $GS_{\bf u,v}^{\bf n,r}(\Delta)$ if for any spline function $p$ belonging to $GS_{\bf u,v}^{\bf n,r}(\Delta)$
\begin{equation}
c_{\xi}=0, \qquad \forall\quad\xi\in {\cal M}\quad \hbox{implies}\quad p\equiv 0, \label{a10ma5}
\end{equation}
where for any $\xi\in {\cal M}$, $c_{\xi}$ is the corresponding B-coefficient of $p$. If there is no smaller set satisfying this property, ${\cal M}$ is called a {\sl minimal determining set}.
\end{definition}
Using simple linear algebra tools, it is easy to verify that, for any determining set ${\cal M}$, $\dim\big(GS_{\bf u,v}^{\bf n,r}(\Delta)\big)\le \# {\cal M}$ and, for any minimal determining set ${\cal M}$, $\dim\big(GS_{\bf u,v}^{\bf n,r}(\Delta)\big)= \# {\cal M}$.

\noindent
Let us denote by $J_{NT}$ and $C$ the sets of vertices which are not T-junctions and the set of composite edges of $\Delta$, respectively. For any $w$ in $J_{NT}$, let $R_w$ be a cell with an edge $e_w$ with one endpoint at $w$, and such that any other edge with an endpoint at $w$ has length at most equal to $e_w$. Moreover, let
\begin{align*}
&{\cal M}_w={\cal D}_{\bf r}^{R_w}(w),\qquad \hbox{for any}\,\,w\in J_{NT},\\ 
&{\cal M}_R=\Big\{\xi_{ij}^R\Big\}_{i=r_1+1,r_2+1}^{n_1-r_1-2,n_2-r_2-2},\qquad \hbox{for any}\,\,R\in \Delta,\\ 
&{\cal M}=\bigcup_{w\in J_{NT}}{\cal M}_w\cup\bigcup_{e\in C}{\cal M}_e\cup \bigcup_{R\in \Delta}{\cal M}_R, 
\end{align*} 
where ${\cal M}_e$, $e\in C$, is defined by \eqref{a8ma1}.

\smallskip
\begin{lemma}\label{a10ma4}
The subset of domain points ${\cal M}\subset {\cal D}_{{\bf n},\Delta}$ is a determining set for $GS_{\bf u,v}^{\bf n,r}(\Delta)$.
\end{lemma}

\smallskip
\noindent
{\bf Proof.} In order to prove the lemma is sufficient to show that \eqref{a10ma5} holds. Let $p\in GS_{\bf u,v}^{\bf n,r}(\Delta)$ with $c_{\xi}=0$ for any $\xi \in {\cal M}$. First of all, by hypothesis, for any $w \in J_{NT}$ $c_\xi=0$ for all $\xi\in {\cal M}_w={\cal D}_{\bf r}^{R_w}(w)$, which implies, by Lemma \ref{a2ma5}, that $c_\xi=0$ for all $\xi \in {\cal D}_{\bf r}(w)$. As a consequence, we can state that for any composite edge $e$ having both the endpoints in $J_{NT}$ we have $c_\xi=0$ for all $\xi \in \tilde {\cal M}_e$, since by hypothesis $c_\xi=0$ for all $\xi \in {\cal M}_e$. Therefore, by using Lemma \ref{a8ma2}, we also obtain that $c_\xi=0$ for all $\xi$ such that $d(\xi,e)\le r_e$. We will refer to the already considered vertices and edges as {\sl already determined}. To determine the remaining B-coefficients we now use an iterative procedure consisting of two steps 
\begin{enumerate}
\item for each T-junction $w$ on an already determined composite edge, we can use Lemma \ref{a2ma5} to show that $c_\xi=0$ for all $\xi\in {\cal D}_{\bf r}(w)$;
\item for each composite edge $e$ with already determined endpoints $w_1$ and $w_2$, since $c_\xi=0$ for all $\xi\in {\cal M}_e$ by hypothesis and $c_\xi=0$ for all $\xi\in {\cal D}_{\bf r}(w_1)\cup {\cal D}_{\bf r}(w_2)$, because these two vertices are already determined, we can use Lemma \ref{a8ma2} to show that $c_\xi=0$ for all $\xi$ such that $d(\xi,e)\le r_e$.
\end{enumerate}
Since we consider only T-meshes without cycles, this two steps can be repeated until all the vertices and all the edges are already determined. Then, at this point, we can state that all the B-coefficients corresponding to domain points within a distance $r_e$ from any edge $e$ are determined and are zero. The remaining coefficients are zeros as well, since they correspond to domain points whose distance from any edge $e$ is greater than $r_e$, which are exactly the domain points contained in $\cup_{R\in \Delta}{\cal M}_R$. \hfill$\square$

\medskip
We are now able to provide the construction of a basis for the space and a dimension formula for the generalized spline spaces $GS_{\bf u,v}^{\bf n,r}(\Delta)$. The logical scheme to get such results is similar  to the one used in \cite{schum12} for the polynomial spline case. Therefore we don't report the proofs.

\medskip
\begin{lemma} \label{a10ma6}
For each $\xi\in{\cal M}$, there is a unique $\psi_\xi\in GS_{\bf u,v}^{\bf n,r}(\Delta)$ such that
\begin{equation*}
\gamma_\eta\psi_\xi=\delta_{\xi,\eta}, \qquad \eta\in{\cal M}, 
\end{equation*}
where $\delta_{\xi,\eta}$ is the Kronecker delta and, for any $\eta\in {\cal D}_{{\bf n},\Delta}$, $\gamma_\eta: GS_{\bf u,v}^{\bf n,r}(\Delta)\rightarrow \RR$ is the functional defined by
\begin{equation}
\gamma_\eta p=c_\eta, \qquad \hbox{with}\,\, c_\eta\,\, \hbox{B-coefficient of $p$ associated to $\xi$}, \quad p\in GS_{\bf u,v}^{\bf n,r}(\Delta). \label{a10ma8}
\end{equation}
\end{lemma}

\medskip
\begin{theorem}\label{a10ma9}
The set ${\cal M}$ is a minimal determining set for $GS_{\bf u,v}^{\bf n,r}(\Delta)$, the set $\{\psi_\xi\}_{\xi\in {\cal M}}$ is a basis for $GS_{\bf u,v}^{\bf n,r}(\Delta)$, and
\begin{align*}
\dim \big( GS_{\bf u,v}^{\bf n,r}(\Delta)\big)=&(r_1+1)(r_2+1)J_{NT}+(r_2+1)(n_1-2r_1-2)E_{hor}\notag\\
&+(r_1+1)(n_2-2r_2-2)E_{ver}+(n_1-2r_1-2)(n_2-2r_2-2)N, 
\end{align*}
where
\begin{align*}
&J_{NT}=\hbox{number of vertices of $\Delta$ which are not T-junctions},\\
&E_{hor}=\hbox{number of horizontal composite edges of $\Delta$},\\
&E_{ver}=\hbox{number of vertical composite edges of $\Delta$},\\
&N=\hbox{number of cells of $\Delta$}.
\end{align*}
\end{theorem}

\medskip
\begin{corollary}\label{a10ma11}
Let ${\bf n}=(n,n)$ and ${\bf r}=(r,r)$ with $n-1\ge 2r+1$; then
\begin{equation*}
\dim \big( GS_{\bf u,v}^{\bf n,r}(\Delta)\big)=(r+1)^2J_{NT}+(r+1)(n-2r-2)(E_{hor}+E_{ver})+(n-2r-2)^2N. 
\end{equation*}
\end{corollary}

\section{Noteworthy cases of generalized spline spaces} \label{sec:400}
\medskip

\noindent Given a T-mesh $\Delta$, it is clear, from \eqref{a2ma1} and \eqref{a2ma2}, that the corresponding generalized spline space $GS_{\bf u,v}^{\bf n,r}(\Delta)$ depends on the choice of the functions $u_1$, $v_1$ and $u_2,v_2$. Some noteworthy choices are, for example, the trigonometric functions, that is, $ u_1(s) = \cos(s)$, $v_1(s) = \sin(s)$, $u_2(t) = \cos(t)$, $v_2(t) = \sin(t) $, or the hyperbolic functions, that is, $ u_1(s) = \cosh(s)$, $v_1(s) = \sinh(s)$, $u_2(t) = \cosh(t)$, $v_2(t) = \sinh(t) $, because of their properties to exactly reproduce certain shapes (conic sections, helices, cycloids, catenaries; see also \cite{manni11}). \\

In particular in isogeometric analysis (see, e.g., \cite{cottrell} and \cite{bazilevs10}), it is of a certain interest to choose $u_1$, $v_1$ and $u_2$, $v_2$ such that the space of generalized splines over the T-mesh has a nice behaviour with respect to the fundamental derivation and integration operators, that is, such that the derivatives and integrals belong to spaces of the same type. Such a feature can be hardly formalized exactly, but in our case we can obtain some favourable cases by examining which spaces $GS_{\bf u,v}^{\bf n,r}(\Delta)$ satisfy the following conditions 
\begin{align}
                         &\psi(s,{t}) \in GS_{\bf u,v}^{\bf n,r}(\Delta) \label{s2_deriv}
\quad\Longrightarrow \quad D_s  \psi(s,{t}) \in GS_{\bf u,v}^{\bf \tilde n_s,\tilde r_s}(\Delta) \\
                         &\psi(s,{t}) \in GS_{\bf u,v}^{\bf n,r}(\Delta) \label{s2_integr}
\quad\Longrightarrow \quad \int \psi(s,{t}) ds \in GS_{\bf u,v}^{\bf \hat n_s,\hat r_s}(\Delta) \\
                         &\psi({s},t) \in GS_{\bf u,v}^{\bf n,r}(\Delta) \label{t2_deriv}
\quad\Longrightarrow \quad D_t  \psi({s},t) \in GS_{\bf u,v}^{\bf \tilde n_t,\tilde r_t}(\Delta) \\
                         &\psi({s},t) \in GS_{\bf u,v}^{\bf n,r}(\Delta) \label{t2_integr}
\quad\Longrightarrow \quad \int \psi({s},t) dt \in GS_{\bf u,v}^{\bf \hat n_t,\hat r_t}(\Delta),                   
\end{align}
where ${\bf\tilde n_s}=(n_1-1, n_2)$, ${\bf\hat n_s}=(n_1+1, n_2)$, ${\bf\tilde n_t}=(n_1, n_2-1)$, ${\bf\hat n_t}=(n_1, n_2+1)$, and ${\bf\tilde r_s}=(r_1-1, r_2)$, ${\bf\hat r_s}=(r_1+1, r_2)$, ${\bf\tilde r_t}=(r_1, r_2-1)$, ${\bf\hat r_t}=(r_1, r_2+1)$.
Although we are dealing with spline spaces in two dimensions defined over a T-mesh, we can study this property in the univariate case over a single interval, because, in every cell, the function is a product between two functions depending just by one variable. Then, from now on we will consider the space $ \mathcal{P}_{u,v}^{n}([a,b]) $, $n\ge 2$, requiring that
\begin{eqnarray}
                         \psi(s) \in \mathcal{P}_{u,v}^n([a,b]) \label{s1_deriv}
& \Longrightarrow &  \psi'(s) \in \mathcal{P}_{u,v}^{n-1}([a,b]) \\
                         \psi(s) \in \mathcal{P}_{u,v}^n([a,b]) \label{s1_integr}
& \Longrightarrow & \int \psi(s) ds \in \mathcal{P}_{u,v}^{n+1}([a,b])                 
\end{eqnarray}

\medskip
First of all, we note that it is not restrictive assuming that the derivatives of $u$ and $v$ are linear combinations of the two functions themselves. In fact, if we require that the derivative of $u$ belongs to the space $ \mathcal{P}_{u,v}^{n-1}([a,b]) $, with $n\ge 3$, we can write it as $ u'(s) = \alpha_1 \cdot u(s) + \alpha_2 \cdot v(s) + \alpha_3 \cdot P_{n-4}(s) $, where $ P_{n-4}(s) $ is a polynomial of degree at most $n-4$, obtained by differentiating the polynomial parts, of degree at most $n-3$, of $u$ and $v$. However, we can assume without loss of generality that both $u$ and $v$ are free from polynomial parts of degree at most $n-3$, because, in the case they are not, the corresponding space $\mathcal{P}_{u,v}^n([a,b])$ would not change. Therefore, we can suppose $ P_{n-4}(s) \equiv 0$ and, as a consequence, $u'(s) = \alpha_1 \cdot u(s) + \alpha_2 \cdot v(s)$. Similar remarks hold for $v$ and for both the functions when considering integration. \\ 
When differentiating a function $\psi(s) \in \mathcal{P}_{u,v}^n([a,b])$, the derivative of its polynomial part of course belongs to $\mathcal{P}_{u,v}^{n-1}([a,b])$. Therefore, condition \eqref{s1_deriv} reduces to

$$ \left\{ \begin{aligned}
u'(s) & = \alpha_1 \cdot u(s) + \alpha_2 \cdot v(s) \\
v'(s) & = \beta_1  \cdot u(s) + \beta_2  \cdot v(s) 
\end{aligned} \right. $$
Following the classical theory of vector-valued ordinary differential equations, we can consider $ A = \begin{bmatrix} \alpha_1 & \alpha_2 \\ \beta_1 & \beta_2 \end{bmatrix} $ and compute its nonzero eigenvalues $ \lambda_1 $ and $ \lambda_2 $, the associated eigenvectors $ W_1 = \begin{bmatrix} w_{11} \\ w_{21} \end{bmatrix} $, $ W_2 = \begin{bmatrix} w_{12} \\ w_{22} \end{bmatrix} $, and examine the possible cases.

\begin{itemize}
\item \textbf{$A$ has distinct real eigenvalues}. In this case, the solutions are given by 
$$ 
\left\{ \begin{aligned}
u(s) & = c_1 w_{11} e^{\lambda_1 s} + c_2 w_{12} e^{\lambda_2 s} \\
v(s) & = c_1 w_{21} e^{\lambda_1 s} + c_2 w_{22} e^{\lambda_2 s}
\end{aligned} \right., 
$$
where $c_1,c_2\in \RR$. Having assumed the linear independence of $u$ and $v$, necessarily we have $ c_1 \neq 0, c_2 \neq 0 $; without loss of generality for the space spanned by $u$ and $v$, we can set 
$$ 
u(s) = e^{\lambda_1 s}, v(s) = e^{\lambda_2 s}
$$ \\
A notable example arises if we set $ \lambda_1 = c, \lambda_2 = -c $: we obtain the same space generated by $ \cosh(cs)$, $\sinh(cs) $. The corresponding Bernstein-like basis with $c=1$, constructed according to Corollary \ref{p1ma11}, is shown in Figures \ref{fig6}-\ref{fig7} for the cases $n=3,5$.

\begin{figure}[H] 
\begin{minipage}[H]{.43\textwidth}
\centering
\includegraphics[scale=0.42]{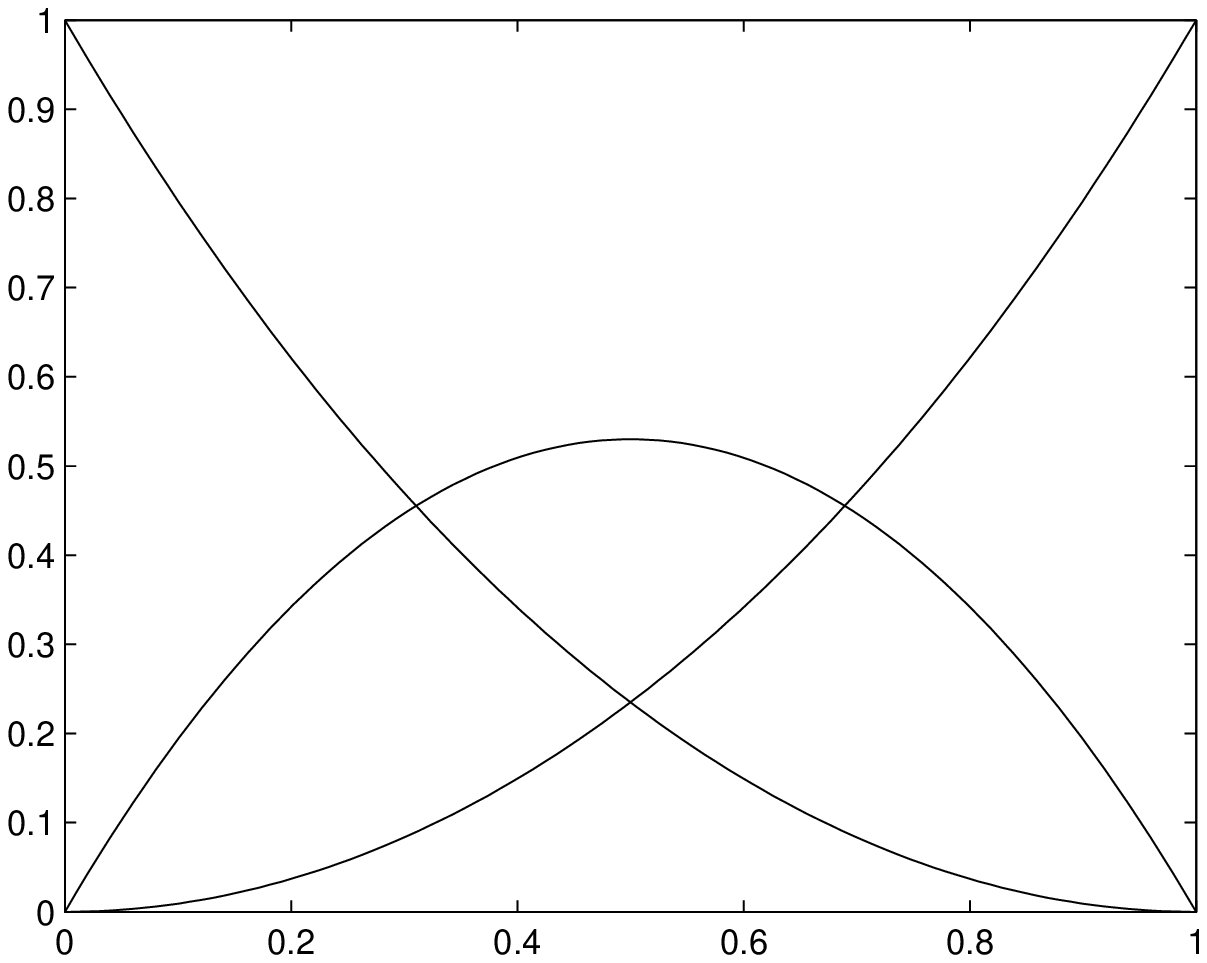} 
\caption{\sl The Bernstein basis of $ \mathcal{P}_{u,v}^{3}([0,1]) $ with $u(s)=\cosh(s)$, $v(s)=\sinh(s)$.}
\label{fig6} 
\end{minipage} 
\hspace{.1\textwidth}
\begin{minipage}[H]{.43\textwidth}
\centering
\includegraphics[scale=0.42]{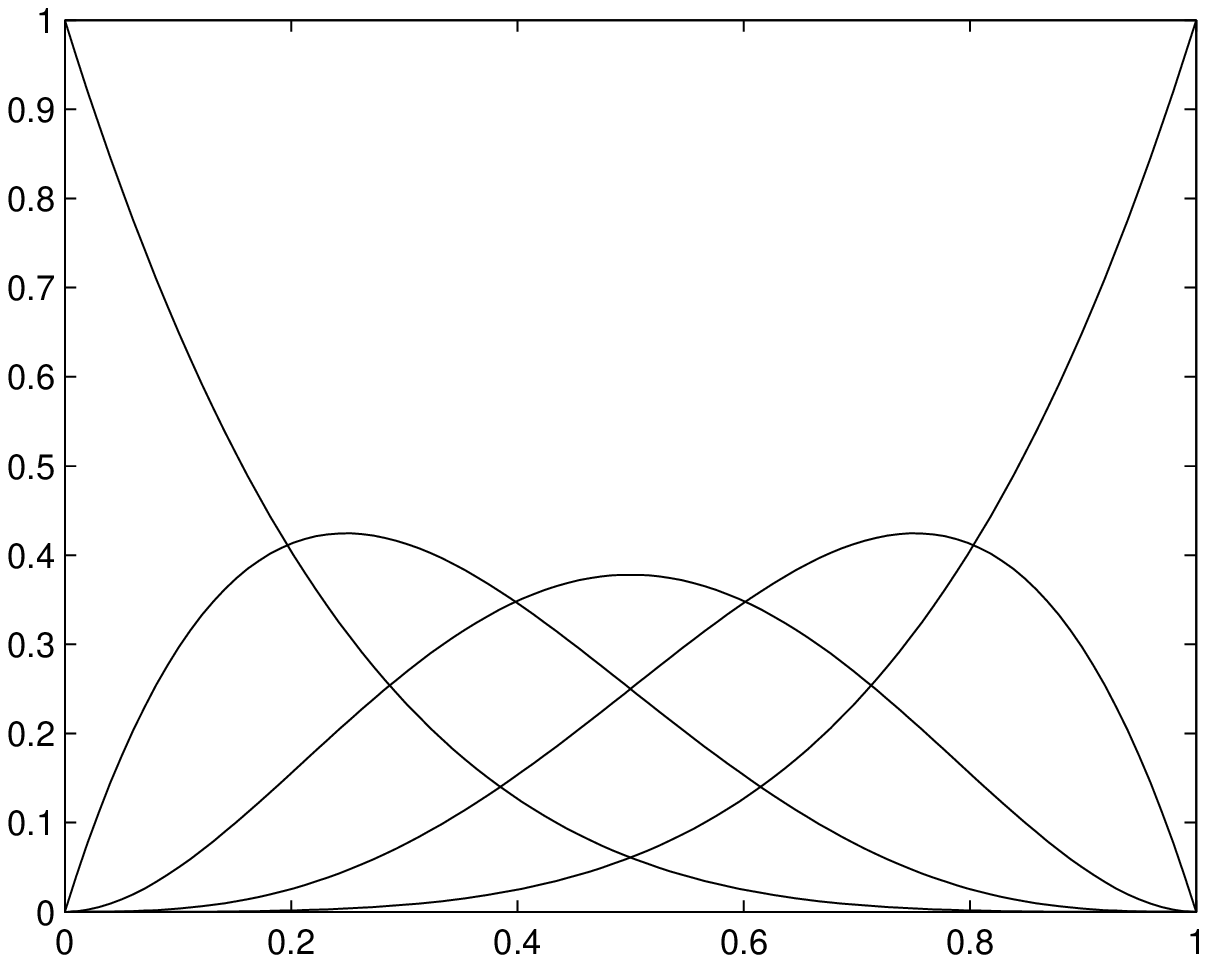}
\caption{\sl The Bernstein basis of $ \mathcal{P}_{u,v}^{5}([0,1]) $ with $u(s)=\cosh(s)$, $v(s)=\sinh(s)$.}
\label{fig7} 
\end{minipage}
\end{figure}

\item\textbf{$A$ has coincident eigenvalues, and it is diagonalizable}. The previous scheme for the solutions can be used, but, due to the fact that $ \lambda_1 = \lambda_2 $, $u$ and $v$ cannot be linearly independent. So this case is not of our interest.

\item\textbf{$A$ has coincident eigenvalues, and it is not diagonalizable}. It corresponds, for our purposes, to $ u(s) = e^{\lambda s}, v(s) = s e^{\lambda s} $.

\item\textbf{$A$ has distinct complex conjugate eigenvalues}. If $ \lambda_1 = \alpha + i \beta, \lambda_2 = \alpha - i \beta $, then 
$$ 
\left\{ \begin{aligned}
u(s) & = c_1 w_{11} e^{(\alpha + i \beta) s} + c_2 w_{12} e^{(\alpha - i \beta) s} \\
v(s) & = c_1 w_{21} e^{(\alpha + i \beta) s} + c_2 w_{22} e^{(\alpha - i \beta) s}
\end{aligned} \right. 
$$
Assuming the linear independence of $u$ and $v$, for our goals it is equivalent considering
$$ 
u(s) = e^{\alpha s} \cos(\beta s), v(s) = e^{\alpha s} \sin(\beta s)
$$
For $ \alpha = 0 $, we obtain the trigonometric case $ u_0(s) = \cos(\beta s), v_0(s) = \sin(\beta s)$. The corresponding Bernstein-like basis with $\beta=1$, constructed according to Corollary \ref{p1ma11}, is shown in Figures \ref{fig8}-\ref{fig9} for the cases $n=3,5$.\\

\begin{figure}[H] 
\begin{minipage}[H]{.43\textwidth}
\centering
\includegraphics[scale=0.42]{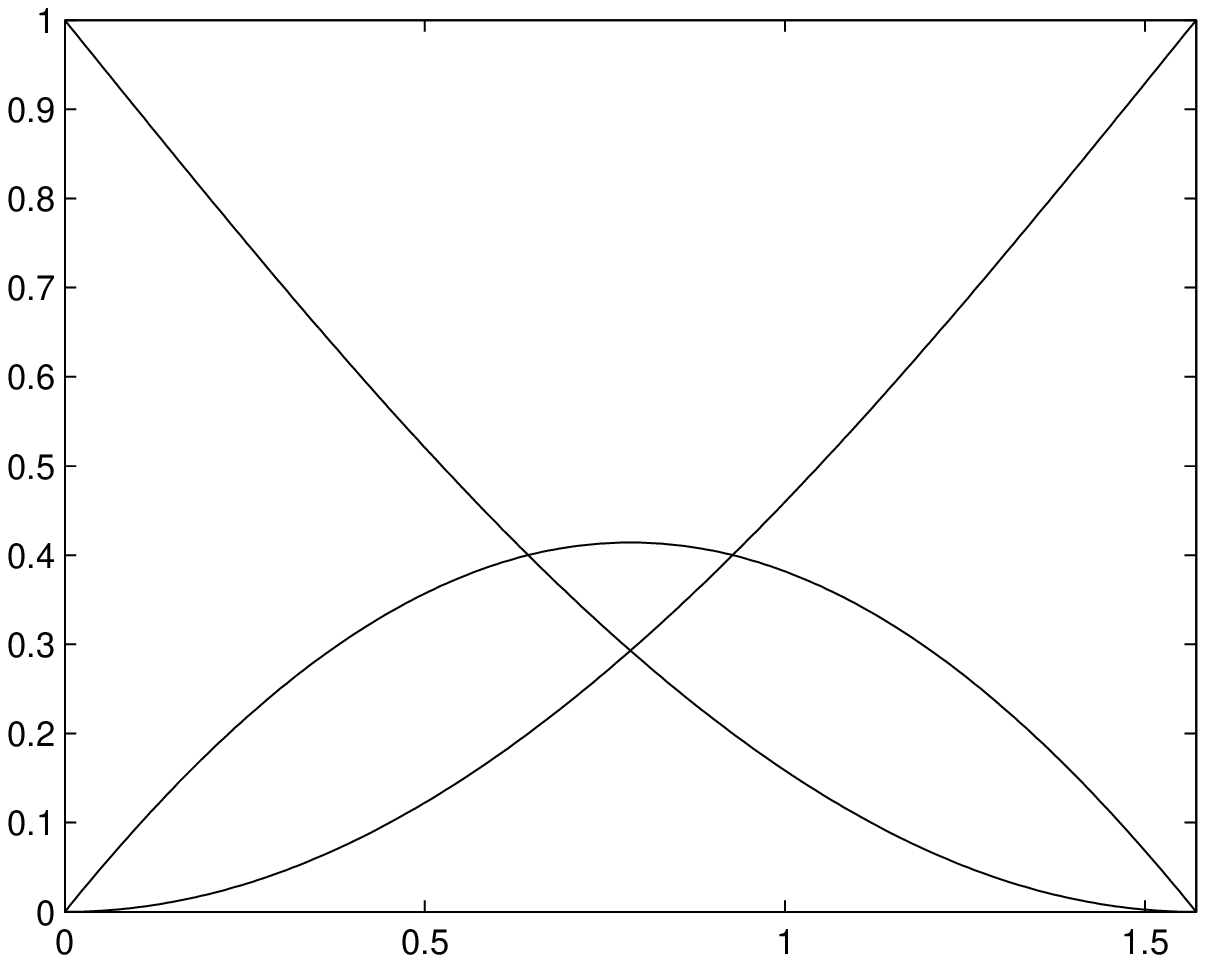}
\caption{\sl The Bernstein basis of $ \mathcal{P}_{u,v}^{3}([0,\pi/2]) $ with $u(s)=\cos(s)$, $v(s)=\sin(s)$.} \label{fig8} 
\end{minipage} 
\hspace{.1\textwidth}
\begin{minipage}[H]{.43\textwidth}
\centering
\includegraphics[scale=0.42]{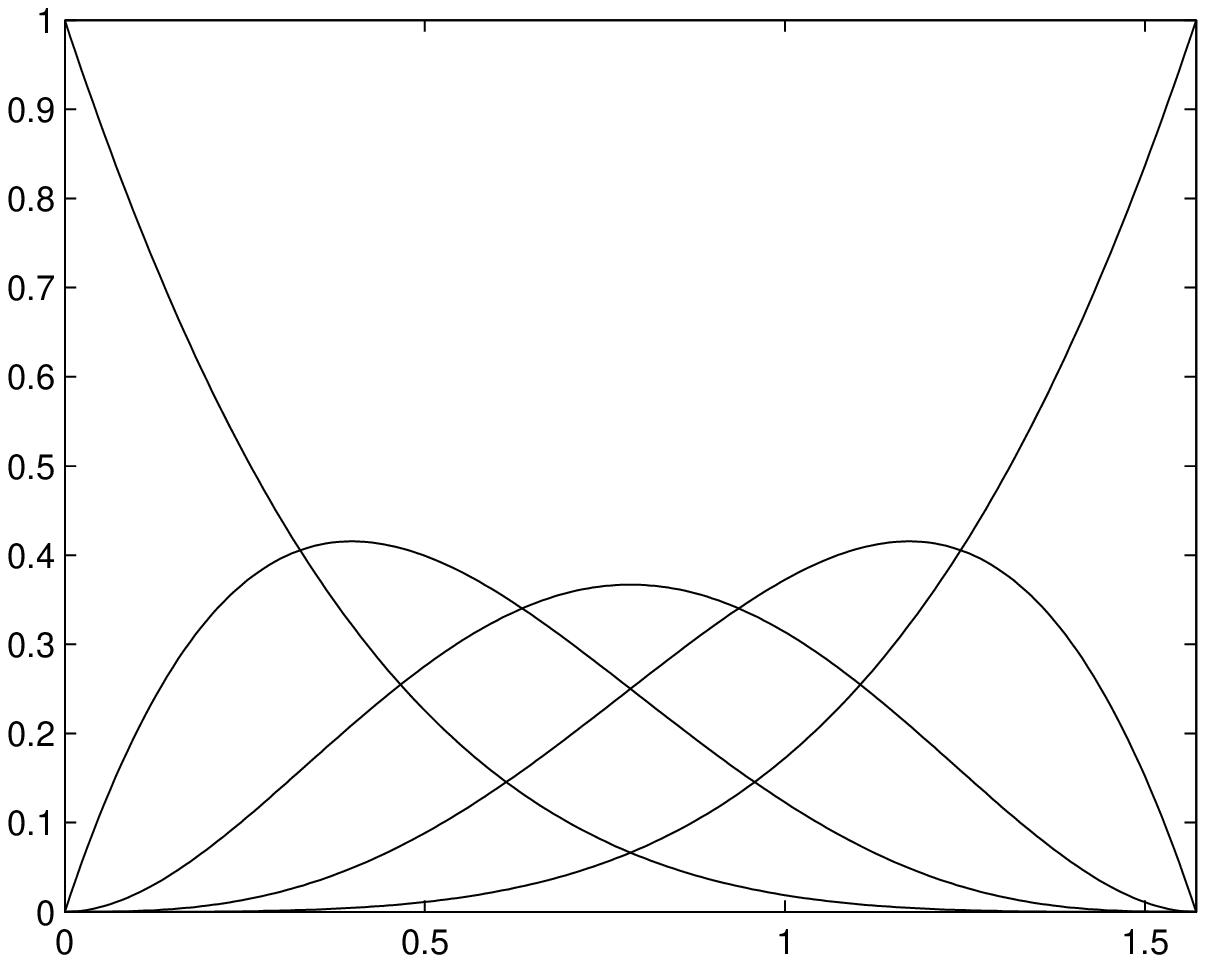}
\caption{\sl The Bernstein basis of $ \mathcal{P}_{u,v}^{5}([0,\pi/2]) $ with $u(s)=\cos(s)$, $v(s)=\sin(s)$.} \label{fig9} 
\end{minipage}
\end{figure}
\end{itemize}

\noindent It can be easily verified that each of the just obtained couples of $u$ and $v$ corresponds to spaces $ \mathcal{P}_{u,v}^n([a,b])$ satisfying \eqref{s1_integr}.
\\ \\

\noindent If \eqref{s1_deriv} and \eqref{s1_integr} hold, \eqref{p1ma2} is verified for sure in some cases and conditionally in other ones, while \eqref{p1ma3} is always verified. By using standard linear algebra means, we obtain, specifically:
\begin{itemize}
\item $ u(s) = e^{\lambda_1 s}, v(s) = e^{\lambda_2 s} $ with $ \lambda_1 \neq \lambda_2 $: \eqref{p1ma2} and \eqref{p1ma3} are both verified for any choice of the parameters and of interval $[a,b]$;
\item $ u(s) = e^{\lambda s}, v(s) = s e^{\lambda s} $: \eqref{p1ma2} and \eqref{p1ma3} are both verified for any choice of the parameters and of interval $[a,b]$;
\item $ u(s) = e^{\alpha s} \cos(\beta s), v(s) = e^{\alpha s} \sin(\beta s) $: \eqref{p1ma2} is verified if $ \beta(b-a) < \pi $, while \eqref{p1ma3} is always verified.
\end{itemize}

\medskip
Finally, we mention another possible choice of $u$ and $v$ (see, e.g., \cite{cost05}) not satisfying the conditions \eqref{s2_deriv} - \eqref{t2_integr}, but anyway having a good behaviour with respect to differentiation and integration 
\begin{equation*} 
u(s) = s^{m_0}, v(s) = (1-s)^{m_1}
\end{equation*} 
where $ m_0 $ and $ m_1 $ are sufficiently large to guarantee the linear independence of the set $\{1,s,...,s^{n-3},s^{m_0},(1-s)^{m_1}\}$. 

\noindent
Also note that all the choices of $u,v$ considered in this section are suitable to construct spaces of generalized splines over a T-mesh \eqref{a2ma1} satisfying the conditions \eqref{p23ot1}-\eqref{p23ot2}.

\section{Approximation power} \label{sec:500}
\medskip

This section is devoted to the study of the approximation properties of the generalized spline spaces over T-meshes. More precisely, we will obtain the approximation order by a logical scheme similar to that one used in \cite{schum12}, but introducing a new suitable quasi-interpolant operator. \\
In fact, the local approximants used in \cite{schum12}, that is, the averaged Taylor expansions, cannot be simply generalized to our non-polynomial case. Moreover, also the results on the approximation power obtained in \cite{cost05} for the univariate case, by using Hermite interpolation in spaces of type ${\cal P}^n_{u,v}([a,b])$, cannot be directly extended to the bivariate case, due to the difficulty to find a suitable differential operator and the corresponding Green's function needed to construct a non-polynomial Taylor expansion. \\
For these reasons, we adopt an alternative approach: we construct a bivariate Hermite interpolant belonging to the spline space, whose existence is rigorously proved by using the assumption \eqref{p1ma2} and \eqref{p1ma3}. This also allows us to obtain an approximation order, which is essentially the same as in polynomial case.

\smallskip
\noindent
Given a function $ f\in C^{\bf n}(\Omega)$ and $(s_0,t_0)\in(a,b)\times (c,d)$, we define the interpolant $Q_L(f;s_0,t_0)(s,t)$ as the function satisfying the two following conditions 
\begin{enumerate}  
\item it belongs to ${\cal P}^{\bf n}_{\bf u,v}([a,b]\times [c,d])$,
\item its polynomial expansion of coordinate bi-degree $(n_1-1,n_2-1)$ coincides with the polynomial expansion of $f$ of the same bi-degree, that is, $Q_L(f;s_0,t_0)(s,t)$ is a Hermite interpolant of coordinate bi-degree $(n_1-1,n_2-1)$. 
\end{enumerate} 
Since $Q_L(f;s_0,t_0)$ is a Hermite interpolant, the Taylor expansion of the difference $f-Q_L(f;s_0,t_0)$ does not contain any term of degree smaller than or equal to $k$, where $ k = \min \{ n_1-1,n_2-1 \} $, and then $ \Vert f-Q_L(f;s_0,t_0) \Vert = O(h^{k+1})$, where $ h = \mbox{diam}([a,b] \times [c,d])$. \\ \\ 
In order to show that $Q_L(f;s_0,t_0)(s,t)$ exists and is unique for any $ f\in C^{\bf n}(\Omega)$ and $(s_0,t_0)\in[a,b]\times [c,d]$, let us write the explicit expressions of a generic element belonging to ${\cal P}^{\bf n}_{\bf u,v}([a,b]\times [c,d])$ 
\begin{align*}
 & \sum_{i=0}^{n_1-3} \sum_{j=0}^{n_2-3} a_{ij} \frac{(s-s_0)^i}{i!} \frac{(t-t_0)^j}{j!} 
           +   \sum_{i=0}^{n_1-3}                    b_i    \frac{(s-s_0)^i}{i!} u_2(t)
           +   \sum_{i=0}^{n_1-3}                    c_i    \frac{(s-s_0)^i}{i!} v_2(t)               \\
          &+                     \sum_{j=0}^{n_2-3} d_j    u_1(s)               \frac{(t-t_0)^j}{j!}
           +                      \sum_{j=0}^{n_2-3} e_j    v_1(s)               \frac{(t-t_0)^j}{j!} \\
          &+  \nu_1 u_1(s) u_2(t) + \nu_2 u_1(s) v_2(t) + \nu_3 v_1(s) u_2(t) + \nu_4 v_1(s) v_2(t)
\end{align*}
and of its Taylor expansion of coordinate bi-degree $(n_1-1,n_2-1)$  
\begin{eqnarray*} \hspace{-10cm} 
&     \sum_{i=0}^{n_1-3} \sum_{j=0}^{n_2-3} a_{ij} \frac{(s-s_0)^i (t-t_0)^j}{i! j!}
              +     \sum_{i=0}^{n_1-3} \sum_{j=0}^{n_2-1} \frac{b_i D_t^j u_2(t_0)}{i! j!} (s-s_0)^i (t-t_0)^j
\\           &+    \sum_{i=0}^{n_1-3} \sum_{j=0}^{n_2-1} \frac{c_i D_t^j v_2(t_0)}{i! j!} (s-s_0)^i (t-t_0)^j
              +     \sum_{i=0}^{n_1-1} \sum_{j=0}^{n_2-3} \frac{d_j D_s^i u_1(s_0)}{i! j!} (s-s_0)^i (t-t_0)^j
\\           &+    \sum_{i=0}^{n_1-1} \sum_{j=0}^{n_2-3} \frac{e_j D_s^i v_1(s_0)}{i! j!} (s-s_0)^i (t-t_0)^j
                \\
								&+\nu_1 \sum_{i=0}^{n_1-1} \sum_{j=0}^{n_2-1} 
                    \frac{D_s^i u_1(s_0) D_t^j u_2(t_0)}{i! j!} (s-s_0)^i (t-t_0)^j \\
		 	 & +    \nu_2 \sum_{i=0}^{n_1-1} \sum_{j=0}^{n_2-1} 
                    \frac{D_s^i u_1(s_0) D_t^j v_2(t_0)}{i! j!} (s-s_0)^i (t-t_0)^j\\
              &+     \nu_3 \sum_{i=0}^{n_1-1} \sum_{j=0}^{n_2-1} 
                    \frac{D_s^i v_1(s_0) D_t^j u_2(t_0)}{i! j!} (s-s_0)^i (t-t_0)^j \\
		 	  &+    \nu_4 \sum_{i=0}^{n_1-1} \sum_{j=0}^{n_2-1} 
                    \frac{D_s^i v_1(s_0) D_t^j v_2(t_0)}{i! j!} (s-s_0)^i (t-t_0)^j .                   
\end{eqnarray*}
\noindent
Then, the condition requiring that $Q_L(f;s_0,t_0)$ is a Hermite interpolant of coordinate bi-degree $(n_1-1,n_2-1)$ corresponds to the following equations 
\begin{align*}
&  a_{ij} + b_i D_t^j u_2(t_0) + c_i D_t^j v_2(t_0) + d_j D_s^i u_1(s_0) + e_j D_s^i v_1(s_0)
  +   \nu_1 D_s^i u_1(s_0) D_t^j u_2(t_0) \\
& +  \nu_2 D_s^i u_1(s_0) D_t^j v_2(t_0) + \nu_3 D_s^i v_1(s_0) D_t^j u_2(t_0) 
  +   \nu_4 D_s^i v_1(s_0) D_t^j v_2(t_0) = D_s^i D_t^j f(s_0,t_0),
\end{align*}
for $ 0 \leq i \leq n_1-3 $, $ 0 \leq j \leq n_2-3 $,
\begin{align*}
&   b_i D_t^j u_2(t_0) + c_i D_t^j v_2(t_0) + \nu_1 D_s^i u_1(s_0) D_t^j u_2(t_0) \\
& +  \nu_2 D_s^i u_1(s_0) D_t^j v_2(t_0) + \nu_3 D_s^i v_1(s_0) D_t^j u_2(t_0) + \nu_4 D_s^i v_1(s_0) D_t^j v_2(t_0) 
  =   D_s^i D_t^j f(s_0,t_0),
\end{align*}
for $ 0 \leq i \leq n_1-3 $, and $ j = n_2-2, n_2-1 $,
\begin{align*}
&    d_j D_s^i u_1(s_0) + e_j D_s^i v_1(s_0) + \nu_1 D_s^i u_1(s_0) D_t^j u_2(t_0) \\
& +  \nu_2 D_s^i u_1(s_0) D_t^j v_2(t_0) + \nu_3 D_s^i v_1(s_0) D_t^j u_2(t_0) + \nu_4 D_s^i v_1(s_0) D_t^j v_2(t_0) 
  =   D_s^i D_t^j f(s_0,t_0),
\end{align*}
for $ i = n_1-2, n_1-1 $, $ 0 \leq j \leq n_2-3 $, and 
\begin{align*}
&    \nu_1 D_s^i u_1(s_0) D_t^j u_2(t_0) + \nu_2 D_s^i u_1(s_0) D_t^j v_2(t_0) \\
& +  \nu_3 D_s^i v_1(s_0) D_t^j u_2(t_0) + \nu_4 D_s^i v_1(s_0) D_t^j v_2(t_0) = D_s^i D_t^j f(s_0,t_0),
\end{align*}
for $i = n_1-2, n_1-1$, $ j = n_2-2, n_2-1 $.
By using a suitable reordering of the unknowns $ a_{ij}, b_i, c_i, d_j, e_j, \nu_k $, we obtain a linear system whose matrix is 
$$ A = \begin{bmatrix}
I & \star & \star & \star \\
0 & A_1 & 0 & \star \\
0 & 0 & A_2 & \star \\
0 & 0 & 0 & A_3
\end{bmatrix} $$ 
where $I$ is the identity matrix of size $ (n_1-2)(n_2-2) \times (n_1-2)(n_2-2) $, $ \star $ stands for blocks of suitable size, $0$ stand for null matrices of suitable size, and
$$ A_1 = \scalemath{0.85}{\begin{bmatrix}
D_t^{n_2-2} u_2(t_0) & 0 & \ldots & 0 & D_t^{n_2-2} v_2(t_0) & 0 & \ldots & 0 \\
D_t^{n_2-1} u_2(t_0) & 0 & \ldots & 0 & D_t^{n_2-1} v_2(t_0) & 0 & \ldots & 0 \\
0 & D_t^{n_2-2} u_2(t_0) & \ldots & 0 & 0 & D_t^{n_2-2} v_2(t_0) & \ldots & 0 \\
0 & D_t^{n_2-1} u_2(t_0) & \ldots & 0 & 0 & D_t^{n_2-1} v_2(t_0) & \ldots & 0 \\
0                  & 0 & \ddots & 0 & 0                  & 0 & \ddots & 0 \\
0 & 0 & \ldots & D_t^{n_2-2} u_2(t_0) & 0 & 0 & \ldots & D_t^{n_2-2} v_2(t_0) \\
0 & 0 & \ldots & D_t^{n_2-1} u_2(t_0) & 0 & 0 & \ldots & D_t^{n_2-1} v_2(t_0) \\
\end{bmatrix}} $$

\medskip
$$ A_2 = \scalemath{0.85}{\begin{bmatrix}
D_s^{n_1-2} u_1(s_0) & 0 & \ldots & 0 & D_s^{n_1-2} v_1(s_0) & 0 & \ldots & 0 \\
D_s^{n_1-1} u_1(s_0) & 0 & \ldots & 0 & D_s^{n_1-1} v_1(s_0) & 0 & \ldots & 0 \\
0 & D_s^{n_1-2} u_1(s_0) & \ldots & 0 & 0 & D_s^{n_1-2} v_1(s_0) & \ldots & 0 \\
0 & D_s^{n_1-1} u_1(s_0) & \ldots & 0 & 0 & D_s^{n_1-1} v_1(s_0) & \ldots & 0 \\
0                  & 0 & \ddots & 0 & 0                  & 0 & \ddots & 0 \\
0 & 0 & \ldots & D_s^{n_1-2} u_1(s_0) & 0 & 0 & \ldots & D_s^{n_1-2} v_1(s_0) \\
0 & 0 & \ldots & D_s^{n_1-1} u_1(s_0) & 0 & 0 & \ldots & D_s^{n_1-1} v_1(s_0) \\
\end{bmatrix} }$$

\medskip
$$ A_3 = \scalemath{0.8}{
\begin{bmatrix}
  D_s^{n_1-2} u_1(s_0) D_t^{n_2-2} u_2(t_0) & D_s^{n_1-2} u_1(s_0) D_t^{n_2-2} v_2(t_0)
& D_s^{n_1-2} v_1(s_0) D_t^{n_2-2} u_2(t_0) & D_s^{n_1-2} v_1(s_0) D_t^{n_2-2} v_2(t_0) \\
  D_s^{n_1-1} u_1(s_0) D_t^{n_2-2} u_2(t_0) & D_s^{n_1-1} u_1(s_0) D_t^{n_2-2} v_2(t_0)
& D_s^{n_1-1} v_1(s_0) D_t^{n_2-2} u_2(t_0) & D_s^{n_1-1} v_1(s_0) D_t^{n_2-2} v_2(t_0) \\
  D_s^{n_1-2} u_1(s_0) D_t^{n_2-1} u_2(t_0) & D_s^{n_1-2} u_1(s_0) D_t^{n_2-1} v_2(t_0)
& D_s^{n_1-2} v_1(s_0) D_t^{n_2-1} u_2(t_0) & D_s^{n_1-2} v_1(s_0) D_t^{n_2-1} v_2(t_0) \\
  D_s^{n_1-1} u_1(s_0) D_t^{n_2-1} u_2(t_0) & D_s^{n_1-1} u_1(s_0) D_t^{n_2-1} v_2(t_0)
& D_s^{n_1-1} v_1(s_0) D_t^{n_2-1} u_2(t_0) & D_s^{n_1-1} v_1(s_0) D_t^{n_2-1} v_2(t_0) \\
\end{bmatrix}} $$
\normalsize The matrix $A_1$ has size $ 2(n_1-2) \times 2(n_1-2) $, $A_2$ has size $ 2(n_2-2) \times 2(n_2-2) $, $A_3$ has size 4. The existence and uniqueness of the interpolation operator $Q_L$ is then equivalent to the non-singularity of this matrix. Since $A$ is an upper triangular block matrix, its non-singularity can be proved by studying $A_1$, $A_2$, $A_3$ ($I$ is obviously non-singular).
\smallskip
\noindent
The matrices $ A_1 $ and $ A_2 $ are not singular, due to their structure, and the fact that \eqref{p1ma2} and \eqref{p1ma3} hold. In fact, the determinants of $A_1$ and $A_2$ are 
\begin{eqnarray*}
\vert \det(A_1)\vert & = &  \vert D_t^{n_2-2} u_2(t_0) D_t^{n_2-1} v_2(t_0) 
                               - D_t^{n_2-1} u_2(t_0) D_t^{n_2-2} v_2(t_0)\vert^{n_1-2} \\
\vert\det(A_2)\vert & = & \vert D_s^{n_1-2} u_1(s_0) D_s^{n_1-1} v_1(s_0) 
                               - D_s^{n_1-1} u_1(s_0) D_s^{n_1-2} v_1(s_0)\vert^{n_2-2}     .                  
\end{eqnarray*}  
Moreover, it can be easily verified that determinant of $A_3$ is $- [\det(D_1)]^2 [\det(D_2)]^2 $, where
\begin{center} \begin{tabular}{ccc}$D_1$ & $:=$ & $ \begin{bmatrix}
D_s^{n_1-2} u_1(s_0) & D_s^{n_1-2} v_1(s_0) \\
D_s^{n_1-1} u_1(s_0) & D_s^{n_1-1} v_1(s_0)
\end{bmatrix} $ \end{tabular} \end{center}
and
\begin{center} \begin{tabular}{ccc}$D_2$ & $:=$ & $ \begin{bmatrix}
D_t^{n_2-2} u_2(t_0) & D_t^{n_2-2} v_2(t_0) \\
D_t^{n_2-1} u_2(t_0) & D_t^{n_2-1} v_2(t_0) 
\end{bmatrix} $ \end{tabular} . \end{center}
If we assume that \eqref{p1ma3} holds, it can be shown (see \cite{cost05}) that $ \det(D_1) \neq 0 $ and $ \det(D_2) \neq 0 $, and so $ \det(A_3) \neq 0 $. \\ In the following, in order to prove the approximation properties, we will also assume that there is a lower bound for the norm of the determinant of $A$ (see also the remarks at the end of the Section). Note that this assumption is guaranteed in the cases where the nonpolynomial functions are $e^{\lambda_1 s}$ and $e^{\lambda_2 s} $, $e^{\lambda s}$ and $ s e^{\lambda s} $, and $e^{\alpha s} \cos(\beta s)$ and $e^{\alpha s} \sin(\beta s) $, examined in Section \ref{sec:400}.

\smallskip
Given a function $f\in C^{\bf n}(\Omega)$, we now define the following quasi-interpolant belonging to the generalized spline space ${GS}^{\bf n,r}_{\bf u,v}(\Delta)$ 
$$ Qf = \sum_{\xi \in \mathcal{M}} \gamma_{\xi}(Q_L(f;s_{\xi},t_{\xi})) \psi_{\xi} $$
where 
\begin{itemize}
\item $ \mathcal{M} $ is the minimal determining set constructed in (34);
\item $ \psi_{\xi} $ are the elements of the basis of the spline space on the T-mesh $\Delta$ associated to $\mathcal{M}$;
\item $ \gamma_{\xi} $ are the linear functionals defined in \eqref{a10ma8} that associate to a spline $p\in {GS}^{\bf n,r}_{\bf u,v}(\Delta)$  the correspondent B-coefficients, needed to express $p$ as a linear combination of the basis $ \psi_{\xi} $: 
\begin{equation*}
p=\sum_{\zeta \in {\cal M}}\gamma_\zeta p \,\psi_\zeta, \qquad \forall p\in {GS}^{\bf n,r}_{\bf u,v}(\Delta) 
\end{equation*}
\item $ (s_{\xi}, t_{\xi}) $ is the center of the biggest circle included in the rectangle $ R_{\xi} $, which is a cell containing $ \xi $. Note that such a point lies in the interior of $ R_{\xi} $, allowing the construction of $ Q_L(f;s_{\xi},t_{\xi}) $. 
\end{itemize}
Note that $Q$ is a linear operator, being the functionals $\gamma_\xi$ linear, and it is a projection onto ${GS}^{\bf n,r}_{\bf u,v}(\Delta)$, that is, $Qp=p$ for every $p\in {GS}^{\bf n,r}_{\bf u,v}(\Delta)$.
In the following, we will use the generalization to our non-polynomial setting of Lemmas 3.1 and 3.2 in \cite{schum12}.
\begin{lemma}\label{C19set1}
Let $p\in {GS}^{\bf n,r}_{\bf u,v}(\Delta)$. Let $R\in \Delta$, and let $p\vert_R=\sum_{\eta \in {\cal D}_{{\bf n}, R}}c_{\eta}^R B_{\eta}^R(s,t)$. We denote by $c$ the vector containing all the coefficients $c^R_\eta$, $\eta \in  {\cal D}_{{\bf n}, R}$. Then, there exists a constant $K_1$, depending only on $n_1$ and $n_2$, such that 
$$ \frac{||c||_{\infty}}{K_1} \leq ||p||_R \leq ||c||_{\infty}, $$
where $ ||c||_{\infty} $ stands for the $\max$-norm of $c$ and $||\cdot||_R $ for the $\sup$-norm of a function restricted to $R$.
\end{lemma}

\smallskip
\noindent
{\bf Proof.}
This is a straightforward generalization of the polynomial case: the upper bound follows from the fact that the basis functions are nonnegative and sum to one, while the lower bound follows from the following argument: the matrix $ M := [B^R_\eta(\zeta)]_{\zeta,\eta\in{\cal D}_{{\bf n}, R}} $ is a non-singular matrix, since the functionals $\{\lambda_\zeta\}_{\zeta\in{\cal D}_{{\bf n}, R}}$, defined by $\lambda_\zeta(f)=f(\zeta)$, $f\in C^0(R)$, are a dual basis of $\{B^R_\eta\}_{\eta\in{\cal D}_{{\bf n}, R}}$. Then $Mc=r$, where $r$ is the vector $\{p(\zeta)\}_{\zeta\in{\cal D}_{{\bf n}, R}}$. As a consequence, we have 
$$ ||c||_{\infty} \leq || M^{-1} r ||_{\infty} \leq || M^{-1} ||_{\infty} ||r||_{\infty} \leq || M^{-1} ||_{\infty} ||p||_R = K_1 ||p||_R .$$
The result is then achieved by setting $ K_1 = || M^{-1} ||_{\infty} $.  \hfill$\square$ 
\smallskip

\begin{lemma}\label{C19set2}
Given a rectangle $R$, let $A_R$ be its area. Then there exists a constant $K_2$, depending only on $n_1$ and $n_2$, such that 
$$ \frac{A_R^{1/q}}{K_2} ||c||_q \leq ||p||_{q,R} \leq A_R^{1/q} ||c||_q, $$
where $ ||c||_q $ stands for the $q$-norm of the vector $c$ and $ ||\cdot||_{q,R} $ for the $q$-norm of a function restricted to $R$.
\end{lemma}

\smallskip
\noindent
{\bf Proof.}
It is sufficient to use equivalence of norms on finite dimensional spaces, considering that both a classical polynomial space and the more general space in which we work have finite dimension. Then, the results is obtained, for any $ 1 \leq q < \infty $, by generalizing Theorem 2.7 in \cite{lai}. \hfill$\square$

\smallskip
To prove the approximation property of the quasi-interpolant, we will need the following result about the minimal determining set and the B-coefficients. 
\smallskip

\begin{definition}\label{aggiuntan}
Let $e$ be a composite edge of $\Delta$, and let $e_1,...,e_m$ be a maximal sequence of composite edges such that for each $i = 1,...,m$, one endpoint of $e_i$ lies in the interior of $e_{i+1}$, where we assume $e_{m+1} = e$. We call $e_1,..., e_m$ a chain ending at $e$. We call $m$ the length of the chain.
\end{definition}

\begin{theorem} \label{sixdotone} 
For every composite edge $e$ consisting of $m$ edge segments $ e_1, \ldots, e_m $ with $ m \geq 1 $, let \linebreak $ \alpha_e := \max \{ |e|/|e_1|, |e|/|e_m| \} $, and let $ \beta_e $ be the length of the longest chain ending on $e$. For each rectangle $R$ in $\Delta$, let $ \kappa_R $ be the ratio of the lengths of its longest and of its shortest edges. Recalling that $ C $ is the set of all composite edges of $\Delta$, we set $ \alpha_{\Delta} := \max_{e \in C} \alpha_e, \beta_{\Delta} := \max_{e \in C} \beta_e, \kappa_{\Delta} := \max_{R \in \Delta} \kappa_R $. Then, for any $ p \in GS^{\bf{n},\bf{r}}_{\bf u,v}(\Delta) $, its associated B-coefficients satisfy 
$$ | c_{\eta} | \leq K_3 \max_{\xi \in \mathcal{M}} |c_{\xi}|,
   \hspace*{1cm} \eta \in {\cal D}_{\mathbf{n}, \Delta} $$
where $K_3$ is a constant depending only on $ \mathbf{n}, \alpha_{\Delta}, \beta_{\Delta}, \kappa_{\Delta} $.
\end{theorem}

\smallskip
\noindent
{\bf Proof.} The proof essentially is a generalization of Theorem 6.1 in \cite{schum12}, since analogous relations between the directional derivatives of a spline belonging to ${GS}^{\bf n,r}_{\bf u,v}(\Delta)$ can be found. \hfill$\square$

\smallskip
\noindent
{\bf Remark.} Let us define, for any cell $R$ in $\Delta$:
\begin{eqnarray*}
\Gamma_R & = & \{ \xi \in \mathcal{M}: \sigma(\psi_{\xi}) \cap R \neq \emptyset \}, \\
\Omega_R & = & \cup_{\xi \in \Gamma_R} \sigma(\psi_{\xi}),
\end{eqnarray*}
where $\sigma(\psi_{\xi})$ denotes the support of $\psi_{\xi}$. Note that, if $ \eta \in {\cal D}_{\mathbf{n},R} $, then 
\begin{equation}\label{lst1} 
\displaystyle | c_{\eta} | \leq K_3 \max_{\xi \in \Gamma_R} |c_{\xi}| 
\end{equation}
since it can be shown that the coefficients corresponding to the domain points $\xi\in \mathcal{M} \backslash \Gamma_R $ do not have influence on the coefficients of ${\cal D}_{\mathbf{n},R}$.

\bigskip
\noindent
Let $\xi\in {\cal M}$ and $F\in C^{\bf n}(\Omega)$. By applying Lemma \ref{C19set2} with $ p = Q_L (F;s_{\xi},t_{\xi}) $, we obtain that $ c = \gamma_{\xi} (Q_L (F;s_{\xi},t_{\xi}) ) $ and
$$ | \gamma_{\xi} (Q_L (F;s_{\xi},t_{\xi}) ) | \leq \frac{K_2}{A_{R_{\xi}}^{1/q}} 
                                                      || Q_L (F;s_{\xi},t_{\xi}) ||_{q,R_{\xi}} .$$ 
If we denote by $T^{(n_1-1,n_2-1)} Q_L(F;s_\xi,t_\xi)$ the Taylor expansion of $Q_L(F;s_\xi,t_\xi)$ at $(s_\xi,t_\xi)$ of bi-degree $(n_1-1,n_2-1)$, for $ 1 \leq q < \infty $, we get 
\begin{eqnarray*}
| \gamma_{\xi} (Q_L (F;s_{\xi},t_{\xi}) ) | & \leq & \frac{K_2}{A_{R_{\xi}}^{1/q}}  
                                            || Q_L (F;s_{\xi},t_{\xi}) ||_{q,R_{\xi}}
                                              \leq   \frac{K_2}{A_{R_{\xi}}^{1/q}} A_{R_{\xi}}^{1/q}
                                            || Q_L (F;s_{\xi},t_{\xi}) ||_{\infty,R_{\xi}} \\
& \leq & K_2 \max_{(s,t) \in R_{\xi}} |T^{(n_1-1,n_2-1)}Q_L(F,s_{\xi},t_{\xi})(s,t)| 
         + O ((\mbox{diam}(R_{\xi}))^{k+1}),
\end{eqnarray*} 
where $ k = \min \{ n_1-1,n_2-1 \} $.
An analogous bound can be obtained for $q=\infty$ by using Lemma \ref{C19set1}. \\ 
For $\eta \in {\cal D}_{{\bf n},R}$, by using Theorem \ref{sixdotone} and \eqref{lst1}, we have 
$$ | \gamma_{\eta} | \leq K_3 \max_{\xi \in \Gamma_R} |\gamma_{\xi}|
                     \leq K_2 K_3 \max_{\xi \in \Gamma_R} \max_{(s,t) \in R_{\xi}}
                          \vert T^{(n_1-1,n_2-1)} Q_L(F;s_\xi,t_\xi)(s,t)\vert
                       +  O(( \max_{\xi \in \Gamma_R} \mbox{diam} (R_{\xi}))^{k+1}) $$ 
That allows us to obtain a limitation for $QF$ 
\begin{eqnarray} 
|| QF ||_{q,R} & \leq & A_R^{1/q} || QF ||_{\infty,R} 
                   =    A_R^{1/q} \left| \left| \sum_{\eta \in {\cal D}_{{\bf n},R}} 
                        \gamma_{\eta} B_{\eta}^R \right| \right|_{\infty,R}\notag \\
               & \leq & A_R^{1/q} K_2 K_3 \max_{\xi \in \Gamma_R} \max_{(s,t) \in R_{\xi}}
                        \vert T^{(n_1-1,n_2-1)}Q_L(F;s_{\xi},t_{\xi})(s,t)\vert\notag\\
                   &+&    O ( ( \max_{\xi \in \Gamma_R} \mbox{diam}(R_{\xi}))^{k+1}) . \label{C19set3} 
\end{eqnarray}
Now, we can finally get an approximation result for the quasi-interpolant $Q$. Given a cell $R_\zeta\in \Delta$, we have 
\begin{eqnarray*} 
|| f - Qf ||_{q,R_{\zeta}} & \leq & || f - Q_L(f;s_{\zeta},t_{\zeta}) ||_{q,R_{\zeta}} +
                                    || Q_L(f;s_{\zeta},t_{\zeta}) - Qf ||_{q,R_{\zeta}} \\
                           &  =   & || f - Q_L(f;s_{\zeta},t_{\zeta}) ||_{q,R_{\zeta}} + 
                                    || Q(f-Q_L(f;s_{\zeta},t_{\zeta})) ||_{q,R_{\zeta}} \\
                          & \leq & O((\mbox{diam}(R_{\zeta}))^{k+1})                        
                              \\&+  &  A_{R_{\zeta}}^{1/q} K_2 K_3
                                   \max_{\xi \in \Gamma_{R_{\zeta}}} \max_{(s,t) \in R_{\xi}}
                                   \vert T^{(n_1-1,n_2-1)}Q_L(f-Q_L(f;s_{\zeta},t_{\zeta});s_{\xi},t_{\xi})(s,t)\vert
               \\[-0.3cm] &  +& \mbox{ } O ( ( \max_{\xi \in \Gamma_{R_{\zeta}}} \mbox{diam}(R_{\xi}))^{k+1}) 
\end{eqnarray*} 
where we used the fact that $Q$ is linear and it is a projection on ${GS}^{\bf n,r}_{\bf u,v}(\Delta)$, and we applied inequality \eqref{C19set3} to $F=f-Q_L(f;s_\zeta,t_\zeta)$. 
Since $\vert D_s^iD_t^j\big(f-Q_L(f;s_{\zeta},t_{\zeta})\big)\vert=O\big(\Vert (s,t) - (s_{\zeta},t_{\zeta}) \Vert^{\max\{0,k+1-i-j\}}\big)$ (for $0\le i\le n_1-1$ and $0\le j\le n_2-1$), and $\Vert (s_{\xi},t_{\xi}) - (s_{\zeta},t_{\zeta}) \Vert\le \mbox{diam}(\Omega_{R_{\zeta}})$, we have 
\begin{align*}
&\vert T^{(n_1-1,n_2-1)}Q_L(f-Q_L(f;s_{\zeta},t_{\zeta});s_{\xi},t_{\xi})\vert\\
&\le \sum_{i=0}^{n_1-1} \sum_{j=0}^{n_2-1} 
\vert  D_s^iD_t^j\big(f-Q_L(f;s_{\zeta},t_{\zeta})\big)\vert_{(s_{\xi},t_{\xi})}  (s-s_{\xi})^i  (t-t_{\xi})^j=O((\mbox{diam}(\Omega_{R_{\zeta}}))^{k+1}).
\end{align*}
Moreover, it can be proved that there exists a constant $K_4$, depending only on $ \alpha_{\Delta}, \beta_{\Delta}, \kappa_{\Delta} $, such that $ \mbox{diam}(\Omega_R) \leq K_4 \mbox{ diam}(R) $ for any $ R \in \Delta $. Then, we get 
\begin{eqnarray*}
|| f - Qf ||_{q,R_{\zeta}} & \leq & O((\mbox{diam}(R_{\zeta}))^{k+1})
                              +     O((\mbox{diam}(\Omega_{R_{\zeta}}))^{k+1})
                              =     O((\mbox{diam}(R_{\zeta}))^{k+1})                            
\end{eqnarray*} 
Then, if we define the \emph{mesh size} of $\Delta$ as $ H = \displaystyle \max_{R \in \Delta} \mbox{diam}(R) $, we obtain  ; 
$$ || f - Qf ||_{q,R_{\zeta}} \leq O(H^{k+1}) $$

\smallskip
\noindent
{\bf Remark.} This approximation order has been obtained assuming that \eqref{p1ma2} and \eqref{p1ma3} are satisfied: in fact this allows to construct the local Hermite interpolants $Q_L (f;s_{\xi},t_{\xi})$. Moreover, we assumed that there is a lower bound for the norm of the determinant of the matrix $A$, which guarantees that the constants involved in the study of the approximation order are bounded. Such conditions are fulfilled in the cases where the nonpolynomial functions are $e^{\lambda_1 s}$ and $e^{\lambda_2 s} $, $e^{\lambda s}$ and $ s e^{\lambda s} $, and $e^{\alpha s} \cos(\beta s)$ and $e^{\alpha s} \sin(\beta s) $, examined in Section \ref{sec:400}. Note that, under analogous assumptions \eqref{p1ma2} and \eqref{p1ma3}, it is possible to construct a univariate version of the Hermite interpolation operator $Q_L$, with approximation order $n$. In other words, for any $f\in C^{n}([a,b])$ and $s_0\in(a,b)$, we can obtain $Q_L(f;s_0)\in {\cal P}^n_{u,v}([a,b])$ satisfying 
\begin{equation}\label{con1}
D_s^kQ_L(f;s_0)(s_0)=D_s^kf(s_0), \qquad k=0,...,n-1,
\end{equation}
and such that $\vert f-Q_L(f;s_0)\vert=O(h^{n})$, where $h=\vert s-s_0\vert$. 

\noindent
This is not in contrast with the approximation results obtained in \cite{cost05} about Hermite interpolation in a space of type ${\cal P}^n_{u,v}([a,b])$. In \cite{cost05} the authors construct, only in the univariate case, the Hermite interpolant in ${\cal P}^n_{u,v}([a,b])$ satisfying conditions \eqref{con1} by using quite different techniques, based on a differential operator (for which ${\cal P}^n_{u,v}([a,b])$ is the nullspace) and on the corresponding Green's function. The approximation results for such interpolant provided in \cite{cost05} are in accordance with ours. In fact, the approximation order obtained under the assumptions \eqref{p1ma2} and \eqref{p1ma3} is $n$ (see Theorem 8 in \cite{cost05}). In one of the examples, where $u(s)=(1-s/h)^\mu$, $v(s) =(s/h)^\mu$ and $[a,b]=[0, h]$ ($h>0$), the authors remark that there is a loss of precision if $\mu > n-1$ when $s_0$ approaches $0$: while the approximation order remains the same, the more $s_0$ gets close to $0$, the more the constant appearing in the error bound grows. However, in our results, this problem never occurs since our assumptions exclude this case.

\section{Conclusions} \label{sec:600}
\medskip

In this paper we provided a deep study of the generalized spline spaces over T-meshes, which generalize the concept of spline space over T-mesh. We showed that, in spite of the different functions locally considered, the overall behaviour of the new spline spaces is analogous to the classical polynomial case. In fact, thanks to the properties of the chosen non-polynomial functions we can use a local Bernstein-B\'ezier representation and generalize the arguments used in \cite{schum12} to the considered non-polynomial case, to get a basis (associated to a minimal determining set) and a dimension formula. Moreover, we obtained the approximation order by considering a quasi-interpolant based on newly defined local Hermite interpolants, which is essentially the same as in the polynomial case. Finally, in Section 4, we also provided a discussion about some noteworthy choices of the non-polynomial functions to be chosen in order to get a good behaviour with respect to differentiation and integration.





\begin{thebibliography}{99}

\bibitem{bazilevs10}
{Y. Bazilevs}, {V.M. Calo}, {J.A. Cottrell}, {J.A. Evans}, {T.J.R. Hughes}, {S. Lipton}, {M.A. Scott} and {T.W. Sederberg}, Isogeometric analysis using T-splines, {Comput. Methods Appl. Mech. Engrg.} {199}, 229-263 (2010).

\bibitem{daveiga}
{L. Beir\~ao da Veiga}, {A. Buffa}, {G. Sangalli} and {R. Vazquez}, Analysis-suitable T-splines of arbitrary degree: definition and properties, {Math. Mod. Meth. Appl. S.} {23}, 1979-2003 (2013).



\bibitem{bracco12b}
{C. Bracco}, {D. Berdisnky}, {D. Cho}, {M. Oh} and {T. Kim}, Trigonometric Generalized T-splines, Comput. Methods Appl. Mech. Engrg. {268}, 540-556 (2014).


\bibitem{cost05}
{P. Costantini}, {T. Lyche}, and {C. Manni}, On a class of weak Tchebycheff systems, {Numer. Math.} {101}, 333-354 (2005).



\bibitem{cottrell}
{J.A. Cottrell}, {T.J.R. Hughes} and {Y. Bazilevs}, Isogeometric Analysis: toward integration of CAD and FEA, John Wiley \& Sons, 2009.

\bibitem{deng06}
{J.S. Deng}, {F.L. Chen}, and {Y.Y. Feng}, Dimensions of spline spaces over T-meshes, {J. Comput. Appl. Math.} {194}, 267-283 (2006).



\bibitem{kvasov99}
{B.I. Kvasov} and {P. Sattayatham}, GB-splines of arbitrary order, {J. Comput. Appl. Math.} {104}, 63-88 (1999).

\bibitem{lai}
{M.J. Lai} and {L.L. Schumaker}, Spline functions on triangulations, Cambridge University Press, 2007.

\bibitem{li06}
{C.J. Li}, {R.H. Wang}, and {F. Zhang}, Improvement on the dimensions of spline spaces over T-meshes, {J. Inform. Comput. Sci.} {3}, 235-244 (2006).

\bibitem{li11}
{X. Li}, {F. Chen}, On the instability in the dimension of spline spaces over T-meshes, {Comput. Aided. Geom. Design} {28}, 420-426 (2011).



\bibitem{manni11}
{C. Manni}, {F. Pelosi} and {M.L. Sampoli}, Generalized B-splines as a tool in isogeometric analysis, {Comput. Methods Appl. Mech. Engrg.} {200}, 867-881 (2011).

\bibitem{manni14}
{C. Manni}, {F. Pelosi} and {H. Speleers}, Local hierarchical $h$-refinements in IgA Based on Generalized B-splines, {Lect. Notes Comput. Sc.} {8177}, 341-363 (2014).



\bibitem{mourrain}
{B. Mourrain}, On the dimension of spline spaces on planar T-meshes, {Math. Comp.} {83}, 847-871 (2014).



\bibitem{schum12}
{L.L. Schumaker} and {L. Wang}, Approximation power of polynomial splines on T-meshes, {Comput. Aided Geom. Design} {29}, 599-612 (2012).

\bibitem{sederberg03}
{T.W. Sederberg}, {J. Zheng}, {A. Bakenov} and {A. Nasri}, T-splines and T-NURCCs, {ACM Trans. Graph.} {22} 3, 477-484 (2003).

\bibitem{sederberg04}
{T.W. Sederberg}, {D.L. Cardon}, {G.T. Finnigan}, {N.S. North}, {J. Zheng} and {T. Lyche}, T-splines simplification and local refinement, {ACM Trans. Graph.} {23} 3, 276-283 (2004).




\end{thebibliography}
\end{document}